\documentclass[11pt]{article}
\usepackage[latin1]{inputenc}    
\usepackage{latexsym} 
\newtheorem{Theorem}{Theorem}[section] 
\newtheorem{lemma}[Theorem]{Lemma} 
\newtheorem{prop}[Theorem]{Proposition} 
\newtheorem{rem}[Theorem]{Remark} 
 
\newtheorem{corol}[Theorem]{Corollary}

\newcommand{\finishproof}{\hfill $\Box$ \vspace{5mm}}

\renewcommand{\theequation}{\thesection .\arabic{equation}}

\begin{document}

\title{KAM Theorem for Gevrey Hamiltonians} 
\author{G. Popov}

\date{} 
\maketitle 
\thispagestyle{empty}

\begin{abstract}

\noindent  
We consider  Gevrey perturbations   $H$ of a  completely integrable 
Gevrey Hamiltonian  $H_0$. Given  a  
Cantor set $\Omega_\kappa$ defined by a Diophantine condition,   
we  find a family   
of KAM invariant tori  of $H$ with frequencies $\omega\in \Omega_\kappa$ 
which is Gevrey smooth in a Whitney sense. 
Moreover, we  obtain a symplectic Gevrey normal form of  
the Hamiltonian in a neighborhood of  the union  $\Lambda$ 
of the invariant  tori.   
This leads to  effective stability of the quasiperiodic motion near 
$\Lambda$.

\vspace{0.3cm}  
\noindent  
\end{abstract}

\setcounter{section}{0}

\section{\it KAM theorem for Gevrey Hamiltonians} 
\setcounter{equation}{0}

Let $D^0$ be a bounded domain in ${\bf R}^n$, and 
${\bf T}^n = {\bf R}^n/2\pi {\bf Z}^n$, $n\ge 2$. We consider a
class of real valued Gevrey Hamiltonians 
in ${\bf T}^n \times D^0$ which 
are small perturbations of a real valued non-degenerate Gevrey
Hamiltonian $H^0(I)$ depending only on the action variables $I\in
D^0$. Our aim is  to obtain a family of KAM (Kolmogorov-Arnold-Moser)
invariant tori $\Lambda_\omega$ 
of $H$ with frequencies $\omega$ in a suitable Cantor
set $\Omega_\kappa$ defined by a Diophantine condition and to prove 
Gevrey regularity for it. It turns out that for each  
$\omega\in \Omega_\kappa$, 
$\Lambda_\omega$ is a   Gevrey smooth embedded torus 
having  the same Gevrey regularity as the Hamiltonian $H$. Moreover,
we shall prove that 
the family $\Lambda_\omega$, $\omega\in \Omega_\kappa$, is Gevrey
smooth with respect to $\omega$ in a Whitney sense, with 
a Gevrey index depending  on the Gevrey class of $H$ and on
the exponent in the Diophantine condition. This naturally involves
anisotropic Gevrey classes. Let $\rho_1,\rho_2 \ge 1$ and $L_1, L_2$
be positive constants.  Given a domain $D \subset {\bf R}^n$, 
we  denote by 
${\cal G}^{\rho_1,\rho_2}_{L_1,L_2}({\bf T}^n \times D)$ 
the set of all $C^\infty$ 
real valued Hamiltonians  $H$ in ${\bf T}^n \times D$  such that 
\begin{equation}
\|H\|_{L_1,L_2} := \sup_{\alpha,\beta\in{\bf N}^n}\, 
\sup_{(\theta, I)\in{\bf T}^n \times D^0  } \, 
\left(|\partial_\theta^\alpha \partial_I^\beta  
H(\theta, I)|\,  L_1^{-|\alpha|} L_2^{-|\beta|}
\alpha !^{-\rho_1} \beta !^{-\rho_2}\right)\ <\ \infty \, ,
                                    \label{a.1}
\end{equation}
where $|\alpha| = \alpha_1 + \cdots + \alpha_n$ and $\alpha ! =
\alpha_1! \cdots \alpha_n!$ for $\alpha=(\alpha_1,\ldots,\alpha_n)
\in {\bf N}^n$. 
In the same way we define 
${\cal G}^{\rho_1,\rho_2}_{L_1,L_2}({\bf T}^n \times \overline D)$, 
where $\overline D$ is the closure of $D$. 
If $\rho_1=\rho_2=\rho$ we write also 
${\cal G}^{\rho}_{L_1,L_2}({\bf T}^n \times D)$, and sometimes we 
do not indicate  the Gevrey   constants $L_1$, $L_2$. 

Let $H^0$ be  a completely 
integrable real valued Gevrey smooth  Hamiltonian 
${\bf T}^n \times D^0 \ni (\theta,I) \to H^0(I)\in {\bf R}$. 
We suppose that $H^0$ is non-degenerate, which means that the map 
$\nabla H^0: D^0 \to \Omega^0$  is a diffeomorphism. Denote by 
$g^0\in C^\infty (\Omega^0)$ the Legendre transform of $H^0$ (then 
$\nabla  g^0: \Omega^0 \to D^0$ is the inverse map to $\nabla  H^0$). 
We suppose
also that there are positive constants  
$\rho >1$,  $A_0 >0$, and $L_0 \le L_2$ such that   
$H^0\in {\cal G}^\rho_{L_0}(D^0)$, 
$g^0 \in {\cal G}^\rho_{L_0}(\Omega^0)$, and 
\begin{equation}
\|H^0\|_{L_0} \, ,\  \|g^0\|_{L_0}\ \le\  A_0 
                                    \label{a.2}
\end{equation}
in the corresponding norms, defined as in (\ref{a.1}). 
In particular, $\Omega^0$ is a bounded
domain. 
Given a subdomain $D$ of $D^0$ we set 
$\Omega := \nabla H^0(D) \subset \Omega^0$. 
Fix $\tau >n-1$ and $\kappa > 0$. We denote by
$\Omega_\kappa$  the set of all 
frequencies $\omega\in \Omega$  having 
distance $\geq\kappa$ to the boundary of $\Omega$ and also 
satisfying the Diophantine condition 
\begin{equation} 
|\langle \omega , k\rangle |\ \geq\ \frac{\kappa}{|k|^\tau}\, ,\;\;\; 
{\rm for\ all}\ 0\neq k\in {\bf Z}^n\, , 
                                  \label{a.3} 
\end{equation} 
where $|k| = |k_1|+\cdots +|k_n|$.

We are going to find a Gevrey family of KAM
invariant tori with frequencies in $\Omega_\kappa$ 
for small perturbations of $H^0$ in ${\cal G}^\rho_{L_1,L_2}$. In what
follows we  fix  the constants $A_0$ and $L_0$, 
and  allow the constants $L_2 \ge L_1\ge 1$
to be arbitrary large. 
This  occurs  in  the case  of the elliptic equilibrium for example 
($L_2 \gg 1$). 
Given $\omega\in \Omega$, we denote by     
${\cal L}_\omega =\langle \omega,\partial_\varphi\rangle = 
\sum_{j=1}^n \omega_j \partial/\partial\varphi_j$ the corresponding 
vectorfield on ${\bf T}^n$. Fix $0<\varsigma\le 1$. 
\begin{Theorem} \label{KAMtheorem}
Let  $H^0 $ be a real valued   
non-degenerate ${\cal G}^\rho$-smooth Hamiltonian , $\rho>1$, depending only
on $I\in D^0$ and  satisfying  (\ref{a.2}). Let $D$ be a subdomain of
$D^0$ with $\overline D \subset D^0$,  
and $\Omega = \nabla H^0 (D)$. Fix $L_2 \ge L_1 \ge 1$ and $\kappa \le 
 L_2^{-1 -\varsigma}$ such that $L_2 \ge L_0$ and $\Omega_\kappa \neq
\emptyset$. 
Then there exists $N=N(n,\rho,\tau) >0$ and  
$\epsilon  >0$ independent  of $\kappa$, $L_1$, $L_2$, 
and of the domain $D\subset D^0$, such that for any 
$H \in {\cal G}^\rho_{L_1,L_2}({\bf T}^n \times \overline D)$ 
 with norm 
$$
\epsilon_{H} := \kappa^{-2} \| H - H^0 \|_{L_1,L_2} \le 
\epsilon L_1^{-N}, 
$$
there exists a  map 
$\overline \Phi: =(\overline U, \overline V): 
{\bf T}^n \times \Omega \to D$ of an 
anisotropic Gevrey class  ${\cal G}^{\rho, \rho'}$, $ \rho'=\rho(\tau
+ 1) +1$,     such that 
\begin{enumerate}
\item[{\rm (i)}]
For each $\omega \in \Omega_\kappa$, 
$\Lambda_\omega := \{(\overline \Phi (\theta,\omega)): \theta\in  {\bf
T}^n\}$ is an embedded  Lagrangian invariant torus of $H$ and  
$X_{H} \circ \overline \Phi(\cdot,\omega)  = 
D\overline\Phi(\cdot,\omega) \cdot {\cal L}_\omega$. 
\item[{\rm (ii)}] There are constants  $A, C >0$, 
independent of  $\kappa$, $L_1$,  $L_2$,  and of  $D$, such
that  
$$
\begin{array}{lrc}
\left|\partial_\theta^{\alpha}\partial_\omega^{\beta}
(\overline U(\theta;\omega) - \theta)\right| + \kappa^{-1}
\left|\partial_\theta^{\alpha}\partial_\omega^{\beta}
(\overline V(\theta;\omega) - \nabla g^0(\omega))\right|\,  \\[0.5cm]
\le\ \displaystyle  A\, 
   C_1^{|\alpha|} 
\left(C_2\kappa^{-1}\right)^{|\beta|} \, 
\alpha !\, ^{\rho } \beta !\, ^{\rho(\tau + 1) +1}\,  
L_1^{N/2}\sqrt{\epsilon_{H}} \, ,
\end{array}
$$
uniformly in  $(\theta,\omega)\in {\bf T}^n\times \Omega$ and for any 
$\alpha,\beta \in {\bf N}^n$, where 
$C_1 = CL_1$ and 
$C_2 =  CL_1^{\tau+1}$.   
\end{enumerate}
\end{Theorem} 
Note that  
$\overline \Phi$ belongs to ${\cal G}^{\rho,\rho'}_{\bar L_1,\bar L_2}
({\bf T}^n \times \overline D)$ with Gevrey constants 
$\bar L_1 = CL_1$ and 
$\bar L_2 =C_2\kappa^{-1} \ge C L_1^{\tau+1} L_2^{1 +\varsigma}$. 
As a consequence, 
we obtain  a
symplectic  normal form of $H$ near the union of the invariant tori. 
 We say that a real valued function 
$\Phi \in C^\infty({\bf R}^n\times D)$ 
is a generating function of an exact symplectic map 
$\chi: {\bf T}^n\times D \to {\bf T}^n\times D$ if 
$\Phi(x,I) - \langle x,I\rangle$ is $2\pi$-periodic with respect to
$x$, 
$|{\rm Id} -
\Phi_I| < 1$,  and 
$$
 \{(\varphi,I;\chi(\varphi,I)):\, (\varphi,I)\in 
{\bf T}^n\times D\} = \{(p(\Phi_I(x,I)),I; p(x), \Phi_x(x,I))
:\, (x,I)\in  {\bf R}^n\times D\},
$$
where $p:{\bf R}^n \to {\bf T}^n$ is the natural projection. 
Fix $\widetilde \kappa = \widetilde \kappa (D)$ so that the Lebesgue
measure of $\Omega_{\widetilde\kappa}$ is positive. 
Define  $\widetilde\Omega_\kappa$, $0< \kappa \le
\widetilde\kappa(D)$, to be the set of  points  
of a positive Lebesgue 
density in $\Omega_\kappa$. In other words, 
$\omega\in\widetilde\Omega_\kappa$ if for  
any neighborhood $U$ of $\omega$ in $\Omega$ the 
Lebesgue measure of $U\cap\Omega_\kappa$ 
is positive. Obviously, $\widetilde\Omega_\kappa$ and 
$\Omega_\kappa$ have the same 
Lebesgue measure. 
\begin{corol} 
Suppose that  
the hypothesis of Theorem 1.1  hold and 
$0< \kappa \le\widetilde \kappa$. Then 
there exists  $N=N(n,\rho,\tau)>0$ and    $\epsilon >0$ 
independent of  $\kappa$, $L_1$,  $L_2$,  and  $D$, such that 
for any 
$H \in {\cal G}^\rho_{L_1,L_2}({\bf T}^n \times \overline D)$ with 
$\varepsilon_H \le \epsilon L_1^{-N-2(\tau + 2)}$ there is 
a ${\cal G}^{\rho'}$-diffeomorphism $\omega:D\rightarrow\Omega$, and an exact 
symplectic transformation  
$\chi\in {\cal G}^{\rho,\rho'}({\bf T}^n\times D,{\bf T}^n\times D) $  
defined by a generating function  $\Phi(x,I)= \langle x,I\rangle +
\phi(x,I)$, 
$\phi\in {\cal G}^{\rho,\rho'}({\bf T}^n\times D)$, such that the
transformed Hamiltonian 
$ \tilde H(\varphi,I) := H(\chi(\varphi,I))$ belongs to 
${\cal G}^{\rho,\rho'}({\bf T}^n\times D)$ and 
for each $I\in \omega^{-1}(\widetilde\Omega_\kappa)$, ${\bf T}^n\times
\{I\}$ is an invariant torus of $\tilde H$. The functions  
$K(I):= \tilde
H(0,I)$ and $R(\varphi,I) : = \tilde H(\varphi,I) - K(I)$ 
satisfy 
$$
\forall \, \alpha
\in {\bf N}^n\ ,\ \forall \, 
(\varphi,I)\in {\bf T}^n\times \omega^{-1}(\widetilde\Omega_\kappa)\
,\quad 
\partial^\alpha \nabla K(I)\, = \,\partial^\alpha \omega(I)\, ,\ 
\partial^\alpha_I  R(\varphi,I)\, =\, 0\, . 
$$ 
Moreover,  there exist $A, C >0$  independent  $\kappa$, $L_1$, 
$L_2$, and of the domain $D$, such that 
$$
\begin{array}{lrc} 
\left| \partial_\varphi^\alpha \partial_I^\beta \phi(\varphi,I) 
\right|\, 
+\, \left|\partial_I^\beta (\omega(I) - \nabla H^0(I))\right|\, +\, 
 \left| \partial_\varphi^\alpha \partial_I^\beta 
(\tilde H(\varphi,I) - H^0(I)) 
\right|  \\ [0.5cm]
\leq\ A \kappa C_1^{|\alpha|} \left(C_2  
\kappa^{-1}\right)^{|\beta|} \, 
\alpha !^{\, \rho} \beta !^{\, \rho'}\, L_1^{N/2} \sqrt{\epsilon_H} \, , 
\end{array}
$$
uniformly with respect to 
$(\varphi,I)\in {\bf T}^n\times D$ and for any $\alpha,\beta \in {\bf
N}^n$, where $C_1 = CL_1$ and 
$C_2 =  CL_1^{\tau+1}$.    
\end{corol}
Denote by $\Omega\ni \omega \mapsto I(\omega)\in D$ the inverse map to
the diffeomorphism $I \mapsto \omega(I)$. Then 
for each $\omega\in \widetilde\Omega_\kappa$,  the restriction
of the Hamiltonian flow of $\tilde H$ to the invariant torus 
${\bf T}^n\times \{I(\omega)\}$ is given by
$(t,\varphi, I) \mapsto (\varphi + t\nabla K(I),I)$, $I = I(\omega)$. 
Set 
$E_\kappa = \omega^{-1}(\widetilde\Omega_\kappa)$.  
Expanding $\partial_\varphi^\alpha \partial_I^\beta R(\varphi, I)$ in
Taylor series at some $I_0\in E_\kappa$ such that 
$|I_0-I| = |E_\kappa -I| = \inf_{I'\in E_\kappa}\, |I'-I|$, we obtain
for any $\alpha,\beta \in {\bf N}^n$ and $m\in {\bf N}$
$$
|\partial_\varphi^\alpha\partial_{I}^{\beta} R(\varphi, I)|\ \leq\ 
\kappa A\,  C_1^{|\alpha|} 
\left(C_2\kappa^{-1}\right)^{|\beta|+m}\, 
\alpha !\,^\rho  \beta !\, ^{\rho'}\, m !\, ^{\rho'-1}, 
\quad  (\varphi,I)\in {\bf T}^n\times D \, ,\ 
I\notin E_\kappa \, ,    
$$ 
where the positive constants $A, C_1,C_2$ are as above. 
Using Stirling's formula 
we minimize the right-hand side  with respect to 
$m\in {\bf N}$ which leads to  
\begin{equation}
|\partial_\varphi^\alpha\partial_{I}^{\beta} R(\varphi, I)|\ \leq\ 
\kappa A\,  C_1^{|\alpha|} 
\left(C_2\kappa^{-1}\right)^{|\beta|}\, 
\alpha !\,^\rho  \beta !\, ^{\rho'}\,  
\exp\left(- ( \kappa C_2^{-1} 
|E_\kappa - I|)^{\, -\frac{1}{\rho(\tau + 1)}}\right)  
                            \label{a.4}
\end{equation}
 for any $\alpha,\beta \in {\bf N}^n$ uniformly with respect to 
$(\varphi,I)\in {\bf T}^n\times D \, ,\ I\notin E_\kappa$, where the
constants $A, C_1, C_2$ 
are as above. 
These inequalities yield   
effective stability of the quasiperiodic motion near 
the invariant  tori as in \cite{kn:Pop}. Effective stability of the
action along all the trajectories for Gevrey smooth Hamiltonians has
been  obtained recently in 
in \cite{kn:M-S}. The importance of the Gevrey category for that  kind
of problems is  indicated by Lochak \cite{kn:Lo}. Integrability
over a Cantor set of tori for 
$C^\infty$ Hamiltonians  is  obtained by P\"{o}schel
\cite{kn:Poe} and
Lazutkin (see \cite{kn:Laz} for references). 

Theorem 1.1 and Corollary 1.2  hold in  the case of a non-degenerate
elliptic equilibrium for Gevrey Hamiltonians as in \cite{kn:Pop}. Indeed, 
let us consider the Birkhoff
normal form of the Hamiltonian, namely, 
$H(\theta,I) = H^0(I) + H^1(\theta,I)$, where 
$H^0(I) = \langle \alpha^0,I\rangle +  \langle Q I,I\rangle $, with 
$\det Q \neq 0$, and $ H^1(\theta,I) = O(|I|^{5/2})$, 
$(\theta, I)$ 
being suitable polar symplectic
coordinates. Here, $ H^1$ is
Gevrey smooth in ${\bf T}^n\times D_a$, $D_a = \{ c_0 a \le I_j \le c_0^{-1}
a:\, j=1,\ldots, n\}$, and $0 < a \le a_0$, where  $0<c_0<1$ is fixed. 
More precisely,  $H^1 \in {\cal G}^\rho_{L_1,L_2}({\bf T}^n\times
D_a)$ with norm  $\|H^1\|_{L_1,L_2} = O(a^{5/2})$, 
where $L_2 = C_0 a^{-1}$, and the positive constants $L_1$ and $C_0$
are fixed. Then  Theorem 1.1 holds choosing  
$\kappa = \delta a^{1+\varsigma}$, $0<\varsigma < 1/4$, $0<\delta \le
1$, and 
for any  $0<a\le a_0 \ll 1$. 

As in 
\cite{kn:Pop} the  symplectic normal form in Corollary 1.2 can be used
to obtain  Gevrey quantum
integrability over the corresponding family of invarint tori  and to 
construct quasimodes with exponentially small
discrepancy in the semi-classical limit for Schr\"{o}dinger type
operators with Gevrey coefficients. Similar results could  be
obtained for more general  classes of non quasi-analytic Hamiltonians
as well. 

The idea of the proof of Theorem 
1.1 is close to that of Theorem 1 in \cite{kn:Pop} (see also
\cite{kn:Kuk}).  It follows from a
KAM theorem for a family of Hamiltonians $P(\theta,I;\omega)$, where 
the frequencies $\omega$ are taken as independent  parameters. Here we
follow closely the exposition of P\"{o}schel \cite{kn:Poe1}. 
First we prove an  approximation lemma for Gevrey
Hamiltonians $P$ with real valued analytic Hamiltonians $P_j$ in
suitable complex domains   in Sect. 3.1. To obtain $P_j$ we first
construct suitable almost analytic extension of $P$ and then we use
Green's formula. In Sect. 3.2 we recall from  P\"{o}schel
\cite{kn:Poe1} the KAM step and in Sect. 3.3 we set the parameters and
make the iterations. Finally, using a Whitney extension theorem due to
Bruna \cite{kn:Bru}, we complete the proof of the theorem. 
In Sect. 3.6 we consider the case of real analytic Hamiltonians and we
improve certain results in \cite{kn:Pop}. 
We prove in
the Appendix an anisotropic version of the implicit function theorem
of Komatsu \cite{kn:Kom} in Gevrey classes.

\section{\it KAM theorem for Gevrey Hamiltonians with parameters} 
\setcounter{equation}{0}

Consider the Hamiltonian $H(\theta,z) = H^0(z) + H^1(\theta,z)$ in 
${\bf T}^n \times D^0$.  
Expanding $H^0(z)$ near given $z_0\in D\subset D^0$, we  write 
$$ 
H^0(z) = H^0(z_0) + \langle \nabla_z H^0 (z_0), 
I\rangle + 
\int_{0}^{1} (1-t) \langle \nabla_z^2 H^0 (z_t)I, 
I\rangle\,  dt, 
$$ 
where $z_t = z_0 + tI$, $z_1=z$, $I$ varies in a small ball 
$B_R(0) = \{|I|<R\}$ in ${\bf R}^n$, and $\nabla_z^2 H^0$ stands for 
the Hessian matrix of $H^0$.  
We put  $\omega = \nabla H^0(z_0)$. Then 
$z_0 = \nabla g ^0(\omega)$, $g^0$ being the Legendre transform of
$H^0$, and we  write 
$$ 
H^0(z) = e(\omega) + \langle \omega, I \rangle + P_{H^0}(I;\omega), 
$$ 
$$ 
H^1(\theta, z) = H^1(\theta,  \nabla g(\omega) + I) = 
P_{H^1}(\theta, I;\omega), 
$$ 
where $e(\omega) = H^0(\nabla g(\omega))$, while 
$P_{H^0}$ stands for the quadratic 
term in $I$ in the expression of $H^0$. We set $P = P_{H^0} +
P_{H^1}$ and consider the  family of 
Hamiltonians 
\begin{equation}
H(\theta,I;\omega) := e(\omega) + \langle \omega, I \rangle
+ P(\theta,I;\omega) 
                                    \label{b.1}
\end{equation} 
in ${\bf T}^n \times B_R(0)$ depending on the frequency $\omega\in
\Omega$. 
From now on, to simplify the notations, we replace 
$C\epsilon_H$, $C A_0$, $CL_1$ and $CL_2$ by 
$\epsilon_H$, $ A_0$, $L_1$ and $L_2$, respectively, whenever $C\ge 1$ 
 depends  only on $L_0$, $\rho$, $\tau$ and $n$. Then 
using (\ref{a.1}), (\ref{a.2}), and Proposition  A.3  we
obtain  
$$
|\partial_\theta^\alpha \partial_I^\beta  \partial_\omega^\gamma 
P(\theta, I;\omega)| \le \left(A_0 R^2 + \kappa^2\epsilon_H\right) 
L_1^{|\alpha|} L_2^{|\beta|+|\gamma|}
(\alpha ! \beta ! \gamma !)^\rho ,    
$$
for any $\alpha,
\beta$ and $\gamma$, and uniformly with respect to $(\theta,I;\omega) \in 
{\bf T}^n \times B_R(0) \times \Omega$.   
Hence, we can suppose that 
$P \in {\cal G}^\rho_{ L_1,L_2,L_2}({\bf T}^n\times B \times
\Omega)$, $B=B_R(0)$,  with norm 
\begin{equation}
\|P\| = \sup  
\left(|\partial_\theta^\alpha \partial_I^\beta \partial_\omega^\gamma 
P(\theta, I;\omega)| L_1^{-|\alpha|} L_2^{-|\beta|-|\gamma|}
(\alpha ! \beta ! \gamma !)^{-\rho}\right) \le   
A_0R^2 + \kappa^2\epsilon_H, 
                                    \label{b.2}
\end{equation}
where the $\sup$ is taken over all multi-indices $\alpha, \beta,
\gamma$ and for all $(\theta,I;\omega) \in 
{\bf T}^n \times B \times \Omega$.  Fix  $0<\varsigma \le  1$.   
We can now formulate our main result in this section. 
\begin{Theorem} \label{KAMtheorem2}
Suppose that $H$ is given by  (\ref{b.1}) where  $P \in {\cal
G}^\rho_{ L_1,L_2,L_2}({\bf T}^n\times \overline 
B_R(0) \times \overline \Omega)$. Fix  
  $\kappa >0$ and $r>0$ such that 
$\kappa, r <  L_2^{-1-\varsigma}$ and $r \le R$.   Then there is
$N=N(n,\rho,\tau)>0$ and 
$\epsilon>0$ independent of $\kappa$, $L_1$, $L_2$, $r$, $R$, and of 
$\Omega \subset \Omega_0$, 
such that if  
\begin{equation}
\| P \| \leq\ \epsilon \kappa r  L_1^{-N}
                                 \label{b.3} 
\end{equation}  
then   
there  exist maps $\phi \in {\cal G}^{\rho '}(\Omega, \Omega)$ and 
$\Phi = (U,V) 
\in {\cal G}^{\rho,\rho '}({\bf T}^n \times\Omega, 
{\bf T}^n \times B_R(0))$, $\rho' = \rho(\tau + 1) + 1$, 
satisfying
\begin{enumerate}
\item[(i)]
For each $\omega \in \Omega_\kappa$, the map 
$\Phi_\omega := \Phi(\cdot,\omega):{\bf T}^n  \rightarrow
 {\bf T}^n\times B_R(0)$ 
is an ${\cal G}^\rho$ embedding and 
$\Lambda_\omega  := \Phi_\omega({\bf T}^n)$ is an embedded Lagrangian 
torus invariant  with respect to the Hamiltonian flow of 
$ H_{\phi(\omega)}(\varphi,I) := H(\varphi, I; \phi(\omega))$. Moreover, 
$X_{H_{\phi(\omega)}} \circ \Phi_\omega = 
D \Phi_\omega \cdot {\cal L}_\omega$  on   ${\bf T}^n$.
\item[(ii)]
There exist $A, C >0$, independent of 
$\kappa$,  $L_1$,  $L_2$, $r$, and of  $\Omega \subset \Omega_0$, such that  
\[
\begin{array}{lrc}
\left|\partial_\theta^{\alpha}\partial_\omega^{\beta}
(U(\theta;\omega) - \theta)\right|\,  + \,
r ^{-1}
\left|\partial_\theta^{\alpha}\partial_\omega^{\beta} V(\theta;\omega)\right|  
+ \kappa^{-1} \left|\partial_\omega^{\beta}
(\phi(\omega) - \omega)\right|
\\    [0.5cm]  
\leq\  A\,  C_1^{|\alpha|}\,  \left(C_2  
\kappa^{-1}\right)^{|\beta|} \, 
\alpha !^{\, \rho} \beta !^{\, \rho'}\,  \displaystyle 
 \frac{\| P \|L_1^{N}} {\kappa r}       
\end{array}
\]
 uniformly in  $(\theta,\omega)\in {\bf T}^n\times \Omega$ and 
for any $\alpha$ and $\beta$, where 
$C_1 = CL_1$ and  $C_2 =  CL_1^{\tau+1}$.    
\end{enumerate}
\end{Theorem} 
{\em Remark}. 
Note that  the constant $e(\omega)$ in  (\ref{b.1})  plays no role in 
Theorem \ref{KAMtheorem2} and from now on we suppose
 $e(\omega) \equiv 0$. 

Theorem \ref{KAMtheorem2} will be proved in the next section. 
Theorem \ref{KAMtheorem} and Corollary 1.2 follow from 
Theorem \ref{KAMtheorem2} and they will be proved in Sect. 4.

\section{\it Proof of Theorem \ref{KAMtheorem2}. } 
\setcounter{equation}{0} 

We divide the proof of Theorem 2.1 in several steps.
First, using Theorem
\ref{Whitney},  we  extend $P$ to a Gevrey function 
$\widetilde P$ of the class 
${\cal G}^\rho_{ \widetilde L_1, \widetilde L_2}
({\bf T}^n\times {\bf R}^{2n})$  
such that $\|\widetilde P\| \le A \|P\|$, $\widetilde L_1 = C L_1$,
and $\widetilde L_2 = C L_2$, where the constants $A,C>0$ are
independent of $P$, $R$ and $\Omega$. 
To simplify the notations we drop $\sim$.
Multiplying $P$  with a suitable cut-off function we
assume  that the support of $P$ with respect to $(I,\omega)$ is
contained in $B_1(0)\times B_{\bar R}(0)$, $\bar R\gg 1$.

\subsection{\it Approximation Lemma for Gevrey functions}
\label{Sec:approximation}
We fix  $0<\varsigma \le 1$ and  choose three strictly 
decreasing sequences of positive numbers
$\{u_j\}_{j=0}^\infty$, $\{v_j\}_{j=0}^\infty$ and
$\{w_j\}_{j=0}^\infty$ tending to $0$ and such  that
\begin{equation}
 \forall j \in {\bf N}\ : \quad   v_jL_2 ,\,   w_jL_2 \le u_j L_1 \le 1\,
,  \quad  v_0 ,\,   w_0\le L_2^{-1-\varsigma}\, . 
                                    \label{c.3}
\end{equation} 
 Consider the complex sets  ${\cal U}_j^m$, $m=1,2$,  in  
${\bf C}^n/2\pi{\bf Z}^n \times  {\bf C}^n \times  {\bf C}^n$ consisting of all
$(\theta, I, \omega)$ with real  parts  ${\rm Re\, } \theta \in {\bf T}^n$, 
${\rm Re\, } I \in B_{2}(0)$ and  ${\rm Re\, } \omega  \in
B_{\bar R+1}(0)$, and such that  $|{\rm Im\, } \theta_k| \le m u_j$, 
 $|{\rm Im\, } I_k| \le m v_j$, $|{\rm Im\, } \omega_k| \le m w_j$, for 
each $1\le k \le n$. Set ${\cal U}_j = {\cal U}_j^1$ and denote by 
$A({\cal U}_j)$ the set of all real-analytic  bounded functions in ${\cal
U}_j$ equipped with the sup-norm $|\cdot|_{{\cal U}_j}$.  

\begin{prop}[Approximation Lemma] \label{approxlemma}
Let $P\in  {\cal G}^\rho_{L_1,L_2}({\bf T}^n 
\times  {\bf R}^{2n})$. Suppose that the  support of $P$ 
with respect to $(I,\omega)$ is 
in $B_1(0)\times B_{\bar R}(0)$, and assume  (\ref{c.3}). 
Then there is a sequence $P_j\in A({\cal U}_j)$, $j\ge 0$, such
that
$$
\begin{array}{rcl}
|P_{j+1} - P_j|_{{\cal U}_{j+1}} &\le&  C_0  \, L_1^n\, 
\exp \left ( - {3\over 4} 
(\rho-1)(2L_1 u_j)^{- {1 \over {\rho -1}}} \right)  \| P \|, 
\\ [0.3cm]
|P_{0}|_{{\cal U}_{0}}  &\le& C_0 \left( 1 + L_1^n
\exp \left ( -{3\over 4} (\rho-1)(2L_1 u_0)^{- {1 \over {\rho -1}}}
\right)\right) \|P\| .
\end{array}
$$
where 
$C_0 = \widetilde C_0(n, \rho,\varsigma )(\bar R^{n} + 1)$.  Moreover, 
$$
\sup |\partial_\theta^\alpha \partial_I^\beta \partial_\omega^\gamma 
(P - P_j)(\theta, I,\omega)| \ \le \ C_0\, L_1^n\,   L_2   \, 
\exp \left ( - {3\over 4} 
(\rho-1)(2L_1 u_j)^{- {1 \over {\rho -1}}} \right) \| P \|, 
$$
in ${\bf T}^n\times B_{1}(0)\times B_{\bar R}(0)$
for  $|\alpha|+ |\beta| + |\gamma| \le 1$.
\end{prop}
{\em Remark}. Instead of $3/4$ we can take above any positive number
less than 1 in order to absorb  certain polynomials of $(L_1
u_j)^{-1}$. Similar estimates can be obtained without the inequalities
$v_0,  w_0 \le L_2^{-1-\varsigma}$. In this case 
$\widetilde C_0= \widetilde C_0(n, \rho)$ but 
 the right hand side of the estimates above should be multiplied 
 by $L_2^{2n}$.  Using the  `standard' proof of the Approximation
Lemma \cite {kn:Ze} one  obtains  for any $\delta >0$ 
an approximation  modulo  
$C(\rho,\delta)\exp \left ( - c(\rho)(L_1 u_j)^{- {1 \over
{\rho}}+\delta} \right) \| P \|$, $C,c>0$. 

\vspace{0.3cm}
\noindent
{\em Proof.} We divide the proof into  two parts. 

\vspace{0.3cm}
\noindent
{\em 1.  Almost analytic extension of P}. 
There is  a constant $C(\rho)\ge 1$, depending only on 
$\rho$, 
such that 
\begin{equation}
 t\in (0,2]\, ,\ m\in {\bf N}\, , \quad 
1 \le m \le t^{\, -{1\over {\rho - 1}}} + 1\,  ,
                                    \label{c.4}
\end{equation}
implies 
\begin{equation}
t^{m}\,  m !^{\, \rho -1}\ \le\ C(\rho) 
m^{(\rho-1)/2}e^{-(\rho-1)m}. 
                                    \label{c.5}
\end{equation}
Indeed, by 
 Stirling's formula, we get
$$ 
t^{m}\,  m !^{\, \rho -1}\ 
\le\  
C_1(\rho) m^{(\rho-1)/2}
e^{-(\rho-1)m} 
\exp\{m[(\rho-1)  \ln m +\ln t ]\}
$$
$$
=\   
C_1(\rho) m^{ (\rho-1)/2}
e^{-(\rho-1)m} 
\exp\left\{(\rho-1) m \ln \left[m\, t^{1 \over {\rho -1}}\right]\right\} .
$$
Moreover, 
$$
m \ln \left[m\, t^{1 \over {\rho -1}}\right] \le m 
\ln \left[1+  t^{1 \over {\rho -1}}\right] \le 
  m\,  t^{1 \over {\rho -1}} \le 1 +   2^{1 \over {\rho
-1}}\, ,
$$ 
which proves (\ref{c.5}).  

We define an almost analytic extensions $F_j$ of $P$ in 
${\cal U}_j^2$ as follows
\begin{equation}
F_j(\theta + i\widetilde \theta,I + i\widetilde I,\omega + i\widetilde \omega)\ 
=\ \sum_{(\alpha,\beta,\gamma)\in
{\cal M}_j}\, 
\partial_\theta^\alpha \partial_I^\beta \partial_\omega^\gamma
P(\theta,I,\omega)\frac{(i\widetilde \theta)^\alpha (i\widetilde I)^\beta 
(i\widetilde \omega)^\gamma}
{\alpha ! \beta ! \gamma !}\ .  
                                    \label{c.6}
\end{equation}
The index set ${\cal M}_j$ consists of all multi-indices 
$\alpha=(\alpha_1,\ldots,\alpha_n)$, 
$\beta=(\beta_1,\ldots,\beta_n)$ and
$\gamma=(\gamma_1,\ldots,\gamma_n)$ 
such that $\alpha_k \le N_1$, 
$\beta_k \le N_2$ and  $\gamma_k \le N_3$, $ k=1,\ldots, n$,   where 
\begin{equation}
N_1 =\left[(2L_1 u_j)^{-{1\over{\rho -1}}}\right] + 1\ ,\ 
N_2 =\left[(2L_2 v_j)^{-{1\over{\rho -1}}}\right] + 1\ , \
N_3 =\left[(2L_2 w_j)^{-{1\over{\rho -1}}}\right] + 1 \  ,
                                    \label{c.7}
\end{equation}
and  $[t]$ stands for the integer part of $t$.  We have 
$$
|F_j|_{{\cal U}_j^2}\ \le \ \|P\| 
\sum_{(\alpha,\beta,\gamma)\in
{\cal M}_j}\, 
(2L_1u_j)^{|\alpha|} (2L_2v_j)^{|\beta|} (2L_2w_j)^{|\gamma|}
(\alpha ! \beta ! \gamma !)^{\rho -1}\ .
$$
For $\alpha_k,\beta_k, \gamma_k \neq 0$ we  
estimate each term 
$$
(2L_1u_j)^{\alpha_k}\alpha_k !^{\, \rho-1}\ ,\  
(2L_2v_j)^{\beta_k}\beta_k !^{\, \rho-1}\ ,\ 
(2L_2w_j)^{\gamma_k}\gamma_k !^{\, \rho -1}\ ,\ k=1,\ldots, n\, ,  
$$ 
by $C(\rho)  m^{(\rho-1)/2}e^{-(\rho-1)m}$, where $m\ge 1$ stands for 
$\alpha_k,\beta_k$, and $\gamma_k$ respectively. To this end, we put 
$t = 2L_1u_j,\, 2L_2v_j,\, 2L_2w_j$, respectively, 
and  we get $t\in
(0,2]$ in view of   (\ref{c.3}).  Now 
(\ref{c.4}) holds because  of   (\ref{c.7}). Then using 
 (\ref{c.5})  we obtain
$$
|F_j|_{{\cal U}_j^2}\ \le \ \|P\|\left(
1 +  C(\rho)\sum_{m=1}^\infty  
m^{(1-\rho)/2}e^{-(\rho-1)m}\right)^{3n} \ =\  
 C_1\|P\| .
$$
Set $z_k =\theta_k + i\widetilde \theta_k$. Then  applying $\bar\partial_{z_k} = 
(\partial/\partial \theta_k + i\partial/\partial\widetilde \theta _k)/2 $ 
to $F_j$  we obtain 
\begin{equation}
2 \bar\partial_{z_k}F_j(\theta + i\widetilde \theta,I + i \widetilde I,
\omega + i\widetilde\omega )\ =\ 
 \sum_{{\scriptstyle{
{(\alpha,\beta,\gamma)\in
{\cal M}_j   }\atop{\alpha_k = N_1 }}}} \, 
\partial_\theta^\alpha \partial_I^\beta \partial_\omega^\gamma
\partial_{\theta_k}
P(\theta,I,\omega)\frac{(i\widetilde \theta)^\alpha (i\widetilde I)^\beta 
(i \widetilde\omega)^\gamma}
{\alpha ! \beta ! \gamma !}\ . 
                                    \label{c.8}
\end{equation}
We  estimate each term in the sum by 
$$
L_1(2L_1u_j)^{|\alpha|} (2L_2v_j)^{|\beta|} (2L_2w_j)^{|\gamma|}
(\alpha ! \beta ! \gamma !)^{\rho -1}(\alpha_k+1)^{\rho} \|P\| 
$$ 
 in   ${\cal U}_j^2$, where $\alpha_k = N_1$.  
Since 
$$
 (2L_1u_j)^{-{1\over{\rho -1}}}  \le \alpha_k = N_1 \le
(2L_1u_j)^{-{1\over{\rho -1}}}  + 1,
$$ 
we obtain  from (\ref{c.5}) (with $t = 2L_1u_j$ and  $m = N_1$)
$$
(2L_1u_j)^{\alpha_k}\,  \alpha_k !^{\, \rho -1} (\alpha_k+1)^{\rho}
\le  C' 
(L_1u_j)^{-{\rho\over{\rho -1}}-{1\over 2}}
\exp\left(-(\rho-1)(2L_1u_j)^{-{1\over{\rho -1}}}\right) .
$$
This implies as above 
$$
\begin{array}{rcl}
|\bar\partial_{z_k}F_j|_{{\cal U}_j^2} \ \le \  
 C'' L_1  (L_1u_j)^{-{\rho\over{\rho -1}}-{1\over 2}}
\exp \left ( - (\rho-1)(2L_1 u_j)^{- {1 \over {\rho -1}}} \right)\,
\|P\| \\  [0.3cm]
|\bar\partial_{z_k}F_j|_{{\cal U}_j^2} \ \le \
C\, L_1\, 
\exp \left ( - \frac{3}{4}\, (\rho-1)(2L_1 u_j)^{- {1 \over {\rho
-1}}} \right)\,  \|P\| \, ,
\end{array}
$$
where $C = C(\rho,n)>0$. 
In the same way, differentiating (\ref{c.8}), we get with $\beta +
\gamma =1$
$$
\begin{array}{rcl}
|\partial_{\theta_p}^\beta\partial_{\widetilde\theta_p}^\gamma
\bar\partial_{z_k}F_j|_{{\cal U}_j^2}\ ,\ 
|\partial_{I_p}^\beta\partial_{\widetilde I_p}^\gamma
\bar\partial_{z_k}F_j|_{{\cal U}_j^2}\ ,\ 
|\partial_{\omega_p}^\beta\partial_{\widetilde \omega_p}^\gamma
\bar\partial_{z_k}F_j|_{{\cal U}_j^2} \\ [0.3cm]
\le C  
L_1 L_2\exp \left ( - \frac{3}{4}\, (\rho-1)(2L_1 u_j)^{- {1 \over {\rho
-1}}} \right)\,  \|P\| \, ,
\end{array}
$$
where $C = C(\rho,n)>0$ (we recall that $L_2 \ge L_1\ge 1$). 
Using (\ref{c.3}), we obtain the same  estimates for 
$\bar\partial_{I_k}F_j$, 
$\bar\partial_{\omega_k}F_j$, and for their  derivatives of order one 
in ${\cal U}_j^2$. 
Indeed, putting $z_k = I_k +
i\widetilde I_k$  we obtain 
$$
\begin{array}{rcl}
|\bar\partial_{z_k}F_j|_{{\cal U}_j^2}  \le  
 C' L_2  (L_2v_j)^{-{\rho\over{\rho -1}}-{1\over 2}}
\exp \left ( - (\rho-1)(2L_2 v_j)^{- {1 \over {\rho -1}}} \right)\,
\|P\|  \\ [0.3cm]
\le C \exp \left ( 
- \frac{3}{4}(\rho-1)(2L_2 v_j)^{- {1 \over {\rho -1}}}
\right)\|P\|  \le  C \exp \left ( 
- \frac{3}{4}(\rho-1)(2L_1 u_j)^{- {1 \over {\rho -1}}}
\right)\|P\| \, ,  
\end{array}
$$
where $C = C(\rho,n,\varsigma)>0$, since $L_2 \le
(L_2v_j)^{-1/\varsigma}$ by (\ref{c.3}). 
We generalize these estimates as follows. Set  
$z=(\theta,I,\omega)\in {\bf C}^{n}/2\pi {\bf Z}^n\times {\bf
C}^{2n}$,  denote by $x_k$ and $y_k$,
respectively, the real and the imaginary part of $z_k$,  $ 1 \le k \le
3n$, and put 
$\bar\partial_{z_k} = (\partial_{x_k} + i \partial_{y_k})/2$.  
Then using (\ref{c.3}), (\ref{c.5}) and (\ref{c.7}), 
we obtain for any   
$\delta = (\delta_1,\ldots, \delta_{3n})\in   {\bf N}^{3n}$ with 
$0\le \delta_k \le  1$ and $|\delta| \ge 1$  the estimate 
$$
|\, \bar\partial_{z}^\delta  F_j\, |_{{\cal U}_j^2} \ \le \ 
C\, L_1^{n}\, 
\exp \left ( - {3\over 4} (\rho-1)(2L_1 u_j)^{- {1 \over {\rho -1}}}
\right)\,   \|P\| ,
$$
where $C = C(n, \rho,\varsigma) >0$. To this end, differentiating 
(\ref{c.6}), we obtain an expression similar to (\ref{c.8}), where
for each $k$ such that  $\delta_k =1$ we have 
$\alpha_k = N_1$ if $1\le k\le n$, $\alpha_k = N_2$ if $n+1\le k\le 2n$,
and $\alpha_k = N_3$ if $2n+1\le k\le 3n$, and then we proceed as above. 
More generally, for any 
$\delta = (\delta_1,\ldots, \delta_{3n})\in   {\bf N}^{3n}$ with 
$0\le \delta_k \le  1$ and
$|\delta| \ge  1$,  and any 
$\beta,\gamma \in   {\bf N}^{3n}$ with $0 \le |\beta| + |\gamma| \le
1$,
  we obtain as above the estimate 
\begin{equation}
| \partial_{x}^\beta 
\partial_{y}^\gamma \bar\partial_{z}^\delta F_j|_{{\cal U}_j^2} \ \le \ 
 C \, L_1^{n}\,  L_2^{|\beta|+ |\gamma|} 
\exp \left ( -  {3\over 4}(\rho-1)(2L_1 u_j)^{- {1 \over {\rho -1}}}
\right)\,   \|P\| ,
                                    \label{c.9}
\end{equation}
where 
$C =  C(n,\rho,\varsigma) >0$. 
Obviously, the same
estimate holds for  $ \bar\partial_{z}^\beta 
\partial_{z}^\gamma \bar\partial_{z}^\delta F_j$ 
if $0 \le |\beta| + |\gamma| \le 1$. 

\vspace{0.5cm}
\noindent
{\em 2.  Construction  of $P_j$}. 
We are going to approximate $F_j$ by analytic  in ${\cal U}_j^2$
functions using   Green's formula
\begin{equation}
\frac{1}{2\pi i} \int_{\partial D}\, \frac{f(\eta)}{\eta -\zeta}\, d\eta  + 
\frac{1}{2\pi i} \int\!\!\!\int_{D}\, 
\frac{\bar\partial_\eta f(\eta)}{\eta -\zeta}\, 
d\eta \wedge d\bar\eta  
= \left\{ \begin{array}{ll} f(\zeta) & \textrm{if $\zeta\in D$}\\
0 & \textrm{if $\zeta\notin \bar D$}, 
\end{array} \right.
                                     \label{c.10}
\end{equation}
where $D \subset \mathbf{C}$ is a bounded domain with a piecewise smooth
boundary $\partial D$ which is positively oriented with respect to
$D$, $\bar D =  D\cup\partial D$, 
and $f\in C^1(\bar D)$.

We denote by $D_k\subset {\bf C}$ the open rectangle 
$\{ |x_k|<a_k,\ |y_k|<b_k \}$, where 
$a_k= \pi$ and   $b_k= 2u_j$ for  $1\le k \le n$;  
$a_k= 2$ and  $b_k= 2v_j$ for  $n+1\le k \le 2n$, and 
$a_k= \bar R +1$ and   $b_k= 2w_j$ for  $2n+1\le k \le 3n$. We denote also 
by  $\partial D_k$ the boundary of $D_k$ which is positively oriented 
with respect to $D_k$ 
and  by $\Gamma$ the union of the oriented segments
$[-\pi - 2i u_j, \pi - 2i u_j] \cup [\pi  + 2i u_j, -\pi + 2i u_j]$. 
Note that $D_k$ and  $\Gamma$ depend on $j$  as well but we omit it. 
Given $\eta\in {\bf C}$,  we consider the $2\pi$-periodic 
meromorphic function 
$$
\zeta \mapsto   K(\eta,\zeta) = \frac{1}{\eta - \zeta} +
K_1(\eta,\zeta), \ K_1(\eta,\zeta) = 
\lim\limits_{N\to +\infty} \, 
\sum\limits_{k=1}^{N} \left(\frac{1}{\eta - \zeta +2\pi k} + 
\frac{1}{\eta - \zeta -2\pi k}\right).
$$
Consider the function
$$
F_{j,1}(z) : = 
\frac{1}{2\pi i} \int_{\Gamma}\, 
F_j(\eta_1,z_2,\ldots,z_{3n})K(\eta_1,z_1)\, d\eta_1 \ , \quad z \in
{\cal U}_j^2 .
$$
It is analytic and $2\pi$-periodic  with respect to $z_1$ in the strip 
$\{|{\rm Im\, }z_1| < 2u_j\}$. Moreover, for $z_1\in D_1$, we have
$$
F_{j,1}(z) = 
\frac{1}{2\pi i} \int_{\partial D_1}\, 
F_j(\eta_1,z_2,\ldots,z_{3n})K(\eta_1,z_1)\, d\eta_1 
$$
since the function under the integral is $2\pi$-periodic with respect to
$\eta_1$, and using (\ref{c.10}) we obtain
$$
F_{j,1}(z) = F_{j}(z) - 
\frac{1}{2\pi i} \int_{D_1}\, 
\bar\partial_{\eta_1} F_j(\eta_1,z_2,\ldots,z_{3n})K(\eta_1,z_1)\, 
d\eta_1 \wedge d\bar\eta_1 .
$$
By continuity last formula remains true for ${\rm Re\, }z_1 = \pm
\pi$. 
Set $F_{j,0}(z): = F_{j}(z)$ and 
${\cal U}_{j,1}^2 : = {\cal U}_{j}^2\cap \{|{\rm Im\, }z_1| \le 
u_j\}$.   We claim that for any multi-index 
$\alpha = (0,\alpha_2,\ldots, \alpha_{3n})\in   {\bf N}^{3n}$ with 
 $0 \le \alpha_m\le 1$, $ 1 \le m \le 3n$, any 
index $k$,  and $\beta, \gamma \in {\bf N}$ such that 
$0 \le \beta + \gamma \le 1$,  
we have 
\begin{equation}
\left |\partial_{z_k}^\beta \bar\partial_{z_k}^\gamma  
\bar\partial_{z}^\alpha 
(F_{j,1} - F_{j,0})\right|_{{\cal U}_{j,1}^2} 
\ \le \ 
 C\, L_1^n\,  L_2^{\beta + \gamma}\, 
\exp \left ( - {3\over 4}(\rho-1)(2L_1 u_j)^{- {1 \over {\rho -1}}}
\right)\,   \|P\| , 
                                    \label{c.11}
\end{equation}
where $C = C(n, \rho,\varsigma) > 0$.
For $k\neq 1$ it directly follows from (\ref{c.9}) differentiating
under the integral. To prove it for
$k=1$, we use the same argument for 
$$
\frac{1}{2\pi i} \int_{D_1}\, \bar\partial_{\eta_1}
F_j(\eta_1,z_2,\ldots,z_{3n})K_1(\eta_1,z_1)\, d\eta_1 \wedge d\bar\eta_1.
$$
On the other hand, 
$$
\frac{1}{2\pi i} \int_{D_1}\, 
\frac{\bar\partial_{\eta_1} F_j(\eta_1,z_2,\ldots,z_{3n})}
{\eta_1 -z_1}\, d\eta_1 \wedge d\bar\eta_1 
$$
$$
= 
- \bar z_1 \bar\partial_{z_1} F_j(z) + \frac{1}{2\pi i} \int_{D_1}\, 
\frac{\bar\partial_{\eta_1} F_j(\eta_1,z_2,\ldots,z_{3n}) -
\bar\partial_{z_1} F_j(z)
}{\eta_1 -z_1}\, d\eta_1 \wedge d\bar\eta_1  
$$
for $z_1 \in D_1$, 
which follows from (\ref{c.10}) applied to $f(z_1) = \bar z_1$. 
Differentiating the last equality 
 and using (\ref{c.9}) we get the
estimate. Moreover, if $|\alpha| \ge 1$, then 
(\ref{c.9}) and (\ref{c.11}) imply 
\begin{equation}
\left |\partial_{z_k}^\beta \bar\partial_{z_k}^\gamma 
\bar\partial_{z}^\alpha
F_{j,1}\right |_{{\cal U}_{j,1}^2} 
\ \le \ 
 C \, L_1^n\,  L_2^{\beta + \gamma}  
\exp \left ( - {3\over 4}
(\rho-1)(2L_1 u_j)^{- {1 \over {\rho -1}}}\right)\, \|P\| . 
                                    \label{c.12}
\end{equation}
We define by recurrence $F_{j,m}(z)$, $2\le m \le n$, 
and we prove that it satisfies (\ref{c.11}) 
in ${\cal U}_{j,m}^2 : = {\cal U}_{j,m-1}^2\cap \{|{\rm Im\, }z_m| \le 
u_j\}$ for $\alpha = (0,\ldots,0,\alpha_{m+1},\ldots,\alpha_{3n})$. 
Moreover, $F_{j,m}(z)$ satisfies (\ref{c.12}) for $|\alpha|\ge 1$. 

For $n<m\le 3n$ we define
$$
F_{j,m}(z) = 
\frac{1}{2\pi i} \int_{\partial D_m}\, 
\frac{F_{j,m-1}(z_1,\ldots, z_{m-1},\eta_m,z_{m+1},\ldots,z_{3n})}
{\eta_m-z_m}\, d\eta_m \ , \quad z \in
{\cal U}_j^2 ,
$$
and set ${\cal U}_{j,m}^2 : = {\cal U}_{j,m-1}^2\cap \{|{\rm Im\, }z_m| \le 
p_m\}$, where $p_m = v_j$ for $n+1 \le m \le 2n$ and
$p_m = w_j$ for $2n+1 \le m \le 3n$. By recurrence with respect to
$m$, we obtain 
(\ref{c.11}) for $F_{j,m}$ 
in ${\cal U}_{j,m}^2$ for any $n<m\le 3n$, 
$\alpha = (0,\ldots,0,\alpha_{m+1},\ldots,\alpha_{3n})$, $0\le
\alpha_q\le 1$  ($\alpha =0$ if $m=3n$), 
for any  index $k$ and 
$\beta, \gamma $ such that $0 \le \beta + \gamma \le 1 $. 
Moreover, 
$F_{j,m}(z)$ satisfies (\ref{c.12}) for $|\alpha|\ge 1$ and  $m<3n$. 
For $2n<m\le 3n$ the constant $C$ should be replaced by 
$C'(1 +{\bar R})^{m-2n}$, where $C'=C'(\rho,n,\varsigma)$. 
The factor $1 +{\bar R}$
comes from the measure of $D_m$, $2n<m\le 3n$.

Set $P_j = F_{j,3n}$. Then for any index $k$
and $\ell =0,1$, we obtain
$$
|\partial_{x_k}^\ell (P_j -F_{j})|_{{\cal U}_{j}} 
\ \le \ 
 C\, L_1^n\,  L_2^\ell
\exp \left ( - {3\over 4} (\rho-1)(2L_1 u_j)^{- {1 \over {\rho -1}}}
\right)  \|P\|  . 
$$
In particular,
$$
\begin{array}{rcl}
|P_{j+1} -P_{j}|_{{\cal U}_{j+1}} \le |P_{j+1} -F_{j+1}| + |P_{j} -F_{j}| +|F_{j+1}
-F_{j}| \\ [0.3cm]
\le C \, L_1^n\, 
\exp \left ( -  {3\over 4} (\rho-1)(2L_1 u_j)^{- {1 \over {\rho -1}}} 
\right)  \|P\| . 
\end{array}
$$
Moreover,  for any index $k$
and $\ell =0,1$,
$$
|\partial_{x_k}^\ell (P_j(x) -P(x))| 
\ \le \  C\, L_1^n\, 
L_2^\ell  
\exp \left ( -{3\over 4} (\rho-1)(2L_1 u_j)^{- {1 \over {\rho -1}}}
\right) \|P\| 
$$
in ${\cal U}_{j}\cap \{ {\rm Im \,} z=0\}  $, since $F_j(x) = P(x)$
for $x$ real. Finally, 
$$
|P_{0}|_{{\cal U}_{0}} \le |F_0| + |P_{0} -F_{0}| \le \  C \left( 1 + L_1^n
\exp \left ( -{3\over 4} (\rho-1)(2L_1 u_0)^{- {1 \over {\rho -1}}}
\right)\right) \|P\| .
$$
This completes the proof of the proposition.
\finishproof

\subsection{\it The KAM step}
\label{Sec:kamstep}

Introduce the complex domains 
$$ 
D_{s,r} = \{\theta\in {\bf C}^n/2\pi {\bf Z}^n:\ |{\rm Im}\,  \theta|<s\} 
          \times \{I\in {\bf C}^n:\ |I|<r\}, 
$$ 
$$ 
O_{h} = \{\omega\in {\bf C}^n:\ |\omega - \Omega_\kappa |<h\}. 
$$ 
The sup-norm of functions in ${\cal V}:=D_{s,r}\times O_h$ will be denoted by 
$|\cdot |_{s,r,h}$. Fix $0< \upsilon <1/6$  
and set $\widetilde \upsilon = 1/2 \,  - 3\upsilon$ (we shall choose later 
$ \upsilon =1/54$ and  $ \widetilde \upsilon = 4/9$ ).  Fix
$0<s,r<1,\ 0<\eta<1/8,\ 0<\sigma < s/5,\ K\geq 1. $
Consider the real valued Hamiltonian $H(\theta,I;\omega) = N(I;\omega) + 
H_1(\theta,I;\omega)$,  
$ N(I;\omega) = e(\omega) + \langle \omega, I\rangle$. 
We shall denote  by `Const.' a positive constant 
depending  only on $n$ and $\tau$ and  by
`const.' if it is $\le 1$. 
We recall from  P\"oschel \cite{kn:Poe1} the following

\begin{prop} \label{kamlemma} Let $H$ be real analytic in ${\cal V}$.  
Suppose that $| H-N |_{s,r,h} \leq \varepsilon$ with 
\begin{itemize}
\item[{\rm (a)}]
$\varepsilon \le {\rm const.}\,  \kappa \eta r \sigma^{\tau + 1},$ 
\item[{\rm (b)}]
$\varepsilon \le {\rm const.}\,\upsilon h r ,$ 
 \item[{\rm (c)}]  $\displaystyle  
 h \leq \frac{\kappa} {2K^{\tau + 1}}$ 
\end{itemize}                            
Then there exists a real analytic transformation 
$$ 
{\cal F}= (\Phi,\phi)\ ,\quad  
\Phi :  D_{s-5\sigma,\eta r}\times O_{h} 
\longrightarrow D_{s,r}\times O_h\, ,\ 
\phi :  O_{\widetilde \upsilon h} \longrightarrow O_h\, ,
$$ 
of the form $\Phi(\theta,I;\omega) = (U(\theta;\omega),
V(\theta,I;\omega))$, with $V$ affine linear with respect to $I$,
where the transformation $\Phi(\cdot ;\omega)$ 
is canonical  for each $\omega$, and 
such that $H\circ {\cal F} = N_+ + P_+$ with 
$N_+(I;\omega)=e_+(\omega) + \langle \omega, I\rangle$, and  
\begin{equation} 
|P_+|_{s-5\sigma,\eta r,\widetilde \upsilon h} 
\le {\rm Const.}\, \left(\frac{\varepsilon ^2}{\kappa r\sigma^{\tau + 1}} 
+ (\eta^2 +K^n e^{-K\sigma})\varepsilon\right) . 
                                   \label{c.13} 
\end{equation} 
Moreover, 
$$ 
|W(\Phi - id )|,\ |W(D\Phi - {\rm Id} )W^{-1}| \le {\rm Const.}\,
 \frac{\varepsilon }{\kappa r\sigma^{\tau + 1}}\, , 
$$ 
$$ 
|\phi - id |,\ \upsilon h|D\phi - {\rm Id} | 
\le {\rm Const.\,} \frac{ \varepsilon }{r}\, ,
$$ 
uniformly on $D_{s-5\sigma,\eta r}\times O_{h}$ 
and  $O_{\widetilde \upsilon h}$, 
respectively, where 
$W = {\rm diag \, }\left(\sigma^{-1}{\rm Id \, }, r^{-1}{\rm Id \,
}\right)$. 
\end{prop} 
\begin{rem}\label{estimates}
Set $\overline W =  
{\rm diag \, }\left(\sigma^{-1}{\rm Id \, }, r^{-1}{\rm Id \, }, 
 h^{-1}{\rm Id \, }\right)$ and suppose that $h \le \kappa
\sigma^{\tau +1}$. Since $1-\widetilde \upsilon \ge 1/3$, 
using the Cauchy estimate with respect to
$\omega$,  we obtain 
$$
|\overline W({\cal F} - id )|,\ |\overline W (D{\cal F}  - {\rm Id} )
\overline W^{-1}| \le  \frac{C \varepsilon }{r h}\, ,
  \ C=C(n,\tau) > 0\, , 
$$
uniformly on $D_{s-5\sigma,\eta r}\times O_{\widetilde \upsilon h}$,  
where $D{\cal F}$ stands for the Jacobian of ${\cal F}$. 
\end{rem}
The proof of Proposition \ref{kamlemma} is given in 
 \cite{kn:Poe1}. The only difference between the statement of
Proposition 1.3 and that of the KAM step in \cite{kn:Poe1} appears in 
 the transformation of the
frequencies ($\widetilde \upsilon = 1/4$ in \cite{kn:Poe1}). 
 To prove the
proposition with $\widetilde \upsilon$ as above we 
use the following analog of Lemma A.3 \cite{kn:Poe1}.
\begin{lemma}
Suppose $f: O_h \rightarrow {\bf C}^n$ is real analytic with bounded 
$|f|_h$. Let $0<\upsilon < 1/6$ and $\widetilde \upsilon = 1/2\,  -
3\upsilon$. If $|f-id|_h \le \upsilon h$, then $f$ has a real analytic
inverse $f: O_{\widetilde \upsilon h} \rightarrow 
O_{2(\upsilon +\widetilde \upsilon) h}$ and 
$$
|\phi -id|_{\widetilde \upsilon h} \ ,\ 
3\upsilon h|D\phi -id|_{\widetilde \upsilon h} \le  |f-id|_h \,  .
$$
\end{lemma}
A sketch of  proof of the Lemma is given in the Appendix. 
We are going to prepare the next iteration. We choose 
a `weighted error' $0<E<1$, fix $0< \widehat\varepsilon \le 1$,  and set 
$$ 
\eta = E^{1/2},\ \varepsilon = 
 \widehat\varepsilon \kappa E r \sigma^{\tau + 1},\ 0<E<1/64 .
$$ 
We define $K$ and $h$ by 
$$ 
K^ne^{-K\sigma} = E,\ h = \frac{\kappa}{2K^{\tau +1}} . 
$$ 
 Setting $x = K\sigma$ we get the equation 
$x^ne^{-x} = E\sigma^n$, 
which has an unique solution with respect to 
$x\in [1,+\infty)$, since 
$0<E<1/64 < 1/e$. Then $K=x\sigma^{-1}>1$.  We set 
$$ 
r_+ = \eta r,\ s_+ = s-5\sigma,\ \sigma_+ = \delta\sigma, 
$$ 
where $0<\delta< 1$. Later we shall choose $\delta =\delta(\rho)$ as a
function of $\rho$ only. 
Now the KAM step gives the estimate 
$$ 
|P_+|_{s_+,r_+,\widetilde \upsilon h} < {\rm Const.}\,  \widehat\varepsilon 
\kappa r\sigma^{\tau + 1}
\left(E^2 
+ (\eta^2 +K^n e^{-K\sigma})E\right) =   {\rm Const.}\, \widehat\varepsilon
 \kappa   r\sigma^{\tau + 1}  E^2  
$$
$$
=  {\rm Const.}\,  \delta(\rho)^{-\tau - 1}\widehat\varepsilon \kappa  
r_+\sigma_+^{\tau + 1}  E^{3/2}. 
$$ 
Hence there is a constant $c_1>1$ depending only on $n$, $\rho$ and $\tau$ 
such that 
$$ 
|P_+|_{s_+,r_+,\widetilde \upsilon h} \leq {1\over 2}
\widehat\varepsilon 
 c_1^{1/2}\kappa r_+\sigma_+^{\tau + 1 }  E^{3/2}. 
$$ 
We fix the weighted error for the iteration by 
$E_+ = c_1^{1/2} E^{3/2}$, set 
$\varepsilon_+ =\widehat\varepsilon \kappa r_+\sigma_+^{\tau + 1}
E_+$,  and then 
define $\eta_+,\ x_+, \ K_+$, and $h_+$ as above. 
Notice that, $c_1E_+ = (c_1 E)^{3/2}$. We require also 
 $c_1E<1$ 
 which leads to an exponentially 
converging scheme. 
Suppose that 
\begin{equation}
h_+ \ \le\  \widetilde \upsilon \, h .
                           \label{c.14}
\end{equation}
Then we obtain 
\begin{equation}
|P_+|_{s_+,r_+,h_+}\ \leq\ {1\over 2} \varepsilon_+. 
                           \label{c.14a}
\end{equation}

\subsection{\it Setting the parameters and iteration}
\label{Sec:iteration}

As in \cite{kn:Poe1} we are going to iterate the KAM step 
infinitely many times choosing appropriately the parameters $0<s,r,\sigma, 
h,\eta<1$ and so on. Our goal is to get a convergent scheme in the Gevrey 
spaces ${\cal G}^{\rho,\rho(\tau +1) +1} $. We are going
to define suitable strictly decreasing sequences of positive numbers 
$\{s_j\}_{j=0}^\infty$, 
$\{r_j\}_{j=0}^\infty$ and  $\{h_j\}_{j=0}^\infty$, tending to $0$,  
and denote 
$$ 
D_j = D_{s_j,r_j}\, ,\ O_j=O_{h_j} \, ,\ {\cal V}_j = D_j\times
O_j .
$$ 
Fix $\delta \in (0,1)$ ($\delta$ will depend only on $\rho$) and set 
$$
s_j = s_0 \delta^j\, ,\ \sigma_j = \sigma_0 \delta^j\, ,\ 
s_0(1-\delta) = 5 \sigma_0 \, .
$$
Obviously, 
$\displaystyle s_{j+1} = s_j - 5\sigma_j$ and 
$\sigma_j = 5^{-1}(1-\delta)s_j$ for $j\ge 0$. 
We set  
$$
u_j = 4s_j = 4s_0 \delta^j \ ,\ v_j = 4r_0\delta^j \ ,\ w_j = 4h_0
\delta^j ,
$$
and denote by ${\cal U}_j$ the corresponding complex sets defined in
Sect. 3.1. We assume for the moment that these sequences verify
(\ref{c.3}). 
Then  applying Proposition  \ref{approxlemma}   we obtain 
\begin{equation}
\begin{array}{rcl}
|P_0|_{{\cal U}_{0}} &\le&   C_0\,  L_1^n\,  \| P \| \,     \label{c.15}\\
|P_{j} - P_{j-1}|_{{\cal U}_{j}}  &\le&  \displaystyle{
 \ C_0\, \,  L_1^n\, \| P \| \,  
e^{ - \widetilde B_0\,  {s_{j}}^{ - {1 \over {\rho -1}}}} = 
\ C_0\, \,  L_1^n\, \| P \| \,  
e^{ - B_0\,  {\sigma_{j}}^{ - {1 \over {\rho -1}}}}\ ,\ 
j\ge 1 \, ,  }
\end{array}
\end{equation}
where the positive constants $\widetilde B_0L_1^{1 \over {\rho -1}}$
and $B_0L_1^{1 \over {\rho -1}}$ depend only on $\rho$ and $\delta$. 
Given $N$ and $a>0$ we  set  
\begin{equation}
\widehat \varepsilon :=  
\| P \|L_1^{N-2} (a  \kappa r )^{-1} \le 1\, , 
                                 \label{c.15a}
\end{equation}
and we  introduce  
$$
\widetilde \varepsilon_{j} =  
\widehat \varepsilon \kappa r_0 \sigma_0^{\tau +1} 
\exp \left( - B_0\,  {\sigma_{j}}^{ - {1 \over {\rho -1}}} \right)\, . 
$$
 We will choose  later $N = N(n,\tau,\rho)$ and  $a>0$ independent of
$\kappa$, $L_1$, $L_2$, and $r$,        so that
$|P_0|_{{\cal U}_0} \le \widetilde \varepsilon_{0}$ and  $
|P_{j} - P_{j-1}|_{{\cal U}_j} \le \widetilde \varepsilon_{j}$ for
$j\ge 1$. Now we  put 
$$
E_j := c_1^{-1}\exp \left( - B {\sigma_{j}}^{ - {1 \over {\rho -1}}} \right)
\quad  {\rm with}\quad B := {B_0\over 2}\left(\delta^{ - {1 \over {\rho
-1}}} - 1\right) ,
$$
where $c_1>1$ is the constant in the KAM step. 
We find  $\delta \in (0,1)$ from the equalities 
$$
\forall j\in {\bf N}\ , \quad  E_{j+1} = c_1^{1/2} E_{j}^{3/2}\, . 
$$
This is equivalent to $\sigma_{j+1} = (2/3)^{\rho-1}\sigma_{j}$, and we get
$\displaystyle \delta = \left(\frac{2}{3}\right)^{\rho-1}$, 
which implies $ \displaystyle  B = B_0/4= 
A_0 {L_1}^{-\frac{1}{\rho-1}}$, where $A_0 = A_0(\rho) >0$ depends
only on $\rho$. 
Now we set $\eta_j = E_j^{1/2}$, 
$\displaystyle r_{j+1} = \eta_j r_{j}$, 
and put
$$
\varepsilon_j =\widehat \varepsilon  \kappa r_j \sigma_j^{\tau +1}  E_j . 
$$
The choice of the `weighted error' $E_j$ above is motivated by the 
inequality $\widetilde \varepsilon_j \le \varepsilon_{j+1}/2$, $j\ge
0$,  which
will be proved in  (\ref{c.18}). 
This inequality will allow us to put 
$P_{j} - P_{j-1}$ in the  error term of the iteration of order $j$.
Next we determine $K_j$ from the equation 
$K_j^ne^{-K_j\sigma_j} = E_j$.  Setting $x_j = K_j\sigma_j$ we obtain
$$
 x_j^n e^{-x_j} = E_j \sigma_j^n = c_1^{-1}\, \sigma_j^n\,  \exp\left(- B
{\sigma_{j}}^{ - {1 \over {\rho -1}}} \right) . 
$$
Consider  the equation
\begin{equation}
x_j - n \ln x_j = B \sigma_{j}^{ - {1 \over {\rho -1}}} - 
n \ln (\sigma_{j}) + \ln c_1 .
                                               \label{c.16}
\end{equation}
We   set
\begin{equation}
\sigma_0\,  =\,   \sigma \,
L_1^{-1}\left(\ln(L_1+e)\right)^{-(\rho-1)}\,  ,\quad  
0< \sigma \le \widetilde \sigma(n,\rho) \ll 1\, . 
                                               \label{c.16a}
\end{equation}
Obviously, $\sigma_0L_1 < \sigma  \le \widetilde \sigma
(n,\rho)  \ll 1$, and  
for  any
$j\in {\bf N}$ we obtain
$$
B  \sigma_{j}^{ - {1 \over {\rho -1}}} - n \ln (\sigma_{j}) + \ln c_1 
\ge B  \sigma_{0}^{ - {1 \over {\rho -1}}}  =  
A_0 (L_1\sigma_{0})^{-\frac{1}{\rho-1}} >  
A_0 \sigma^{-\frac{1}{\rho-1}} \gg1 \, .
$$
Hence,  choosing $0<\sigma  \le \widetilde \sigma (n,\rho) \ll 1$, 
we obtain for each $j\in {\bf N}$ 
an unique solution $x_j = x_j(\sigma)$  of 
(\ref{c.16}) such that 
$$
x_j \ge x_j -n\ln x_j \ge  
B \sigma_{j}^{ - {1 \over {\rho -1}}} 
\ge A_0 \sigma^{ - {1 \over {\rho -1}}}\gg 1\, .
$$    
Then $x_j -n\ln x_j = x_j (1 + o(1))$ as 
$\sigma\searrow 0$. On the other hand,  using again (\ref{c.16a}) we get
$$
\begin{array}{rcl}
x_j -n\ln x_j &\le&   B  \sigma_{j}^{ - {1 \over {\rho -1}}}
\left[ 1 - n A_0^{-1}(L_1\sigma_{j})^{\frac{1}{\rho-1}} \ln(L_1\sigma_{j})
+  n A_0^{-1}(L_1\sigma_{0})^{\frac{1}{\rho-1}} (\ln L_1 +\ln(c_1))
\right]    \\[0.3cm]
&=& B  \sigma_{j}^{ - {1 \over {\rho -1}}} (1 + o(1))\, ,
\end{array}
$$
uniformly with respect to $j\in {\bf N}$.  Hence, 
\begin{equation}
B \sigma_{j}^{ - {1 \over {\rho-1}}} \le x_j \le   
B  \sigma_{j}^{ - {1 \over {\rho -1}}}(1 + o(1))\, ,\  \sigma\searrow 0, 
                                    \label{c.17}
\end{equation}
uniformly with respect to $j\in {\bf N}$. 
We  set 
$ \displaystyle h_j = \kappa\,  2 ^{-1} K_j^{-\tau -1}$ and fix 
$\upsilon = 1/54$. 

We are going to check the hypothesis (a) and (b)
in  Proposition \ref{kamlemma} 
for any $j\ge 0$ ( (c) is fulfilled by definition). 
To prove  $(a)$ we use (\ref{c.15a}) and that $\eta_j^2 = E_j =o(1)$
as $\sigma \searrow 0$. 
Using  (\ref{c.16a}) and (\ref{c.17})  we obtain  
$$
\begin{array}{rcl}
\frac{\varepsilon_j}{r_j h_j}\,  =\,    
2\widehat \varepsilon  E_j {x_j}^{\tau + 1}\,   \le  \, 
2 c_1^{-1}\exp 
\left ( -  B {\sigma_{j}}^{ - {1 \over {\rho -1}}}\right)
\left (B {\sigma_{j}}^{ - {1 \over {\rho -1}}}\right)^{\tau +1}(1 +
o(1)) \\ [0.3cm]
\le c(\rho,\tau)\left ( -  {A_0\over 2} {(L_1\sigma_{j})}^{ - {1 \over
{\rho -1}}}\right)  \le  
c(\rho,\tau)\left ( -  {A_0\over 2} {(\sigma \delta^{j})}^{ - {1 \over
{\rho -1}}}\right).
\end{array}
$$
This  implies 
$\varepsilon_j(r_j h_j)^{-1}  \ll {\rm
const.}\, \upsilon$ for $0<\sigma  \le \widetilde \sigma (n,\rho,\tau)
\ll 1$ which proves (b). In the same way  we obtain 
\begin{equation}
\prod_{j=0}^\infty 
\left( 1 + \frac{C \varepsilon_j}{r_jh_j}\right) \, \le \, 
\exp\left(\sum_{j=0}^\infty 
 \frac{C \varepsilon_j}{r_jh_j}\right) \, \le \, 2
                                    \label{c.17a}
\end{equation}
for  $0<\sigma  \le \widetilde \sigma (n,\rho,\tau) \ll 1$, where
$C=C(n,\tau) >0$ is the constant in Remark 3.3. 
We are going  to check (\ref{c.14}) with $\upsilon =
1/54$. 
Using (\ref{c.17}) we obtain 
$x_{j}/x_{j+1} = 
\left(\sigma_{j+1}/\sigma_{j}\right)^{  {1 \over {\rho -1}}}
(1 + o(1))$. 
This  implies  
$$
\frac{h_{j+1}}{h_{j}} =
\left(\frac{x_{j}}{x_{j+1}}\right)^{\tau+1}
\left(\frac{\sigma_{j+1}}{\sigma_{j}}\right)^{\tau+1} = 
\delta^{(\tau+1){\rho\over{\rho -1}}}
 \left(1 + o(1))\right) 
= \left(\frac{2}{3}\right)^{\rho(\tau + 1)} (1 + o(1))
$$
for  $\sigma \searrow 0$,  uniformly with respect to $j\in {\bf N}$. 
Since $\rho >1$ and $\tau
+1> n \ge 2$,  
we obtain 
$$
\frac{h_{j+1}}{h_{j}} < 
\left(\frac{4}{9}\right)^{\rho} < \frac{4}{9}  
= \widetilde \upsilon , 
$$
for any   $0 < \sigma  \le \widetilde \sigma (n,\rho,\tau) \ll 1$, 
which  proves (\ref{c.14}). 

Using the special choice of $E_j$, we  are going to prove by induction
that
\begin{equation}
\forall\, j\in {\bf N}\ , \quad \widetilde \varepsilon_j\  \le \
{1\over 2} \varepsilon_{j+1} .
                                     \label{c.18}
\end{equation}
To obtain the estimate for $j=0$ we write 
$\varepsilon_{1} =\widehat\varepsilon \kappa  r_{1} \sigma_{1}^{\tau +1} 
 E_{1}
=\widehat\varepsilon  \kappa  r_{0}\sigma_{0}^{\tau
+1} \delta^{\tau +1}c_1^{1/2} E_{0}^{2}$, and we obtain
$$
\widetilde \varepsilon_0/\varepsilon_1 = 
c_1^{-1/2} \left(\frac{3}{2}\right)^{(\tau + 1)(\rho-1)}
\exp \left( - 2B {\sigma_{0}}^{ - {1 \over {\rho -1}}} \right) \le 1/2
$$
for $0 < \sigma  \le  \widetilde \sigma (n,\rho,\tau) \ll 1$, 
since $B_0 = 4B$. 
To prove it for $j +1\ge 1$ we write 
$$
\varepsilon_{j+2} =\widehat\varepsilon \kappa  r_{j+2} \sigma_{j+2}^{\tau +1} 
 E_{j+2}
=\widehat\varepsilon  \kappa  
 r_{j+1}\sigma_{j+1}^{\tau
+1} \delta^{\tau +1}c_1^{1/2} E_{j+1}^{2}
= \varepsilon_{j+1}c_1^{1/2} E_{j+1}\delta^{\tau + 1} .
$$ 
Then for  $j\ge 0$ we obtain
$$ 
\frac{\widetilde\varepsilon_{j+1}}{\varepsilon_{j +2}} \, 
\left(\frac{\widetilde\varepsilon_{j}} {\varepsilon_{j+1}}\right)^{-1} 
= c_1^{-1/2} \left(\frac{3}{2}\right)^{(\tau + 1)(\rho-1)}
\exp \left( -\frac{1}{2} B {\sigma_{j}}^{ - {1 \over {\rho -1}}}
\right)  
\le 
 1
$$
for $0 < \sigma \le  \widetilde \sigma (n,\rho,\tau) \ll 1$
which implies by recurrence (\ref{c.18}). 
From now on we fix $\sigma=  \sigma (n,\rho,\tau) \ll 1$ so that all the
estimates above hold and define  $\sigma_0$ by (\ref{c.16a}). 
Then we set $s_0= 5\sigma_0 (1-\delta)^{-1}$. 
We are going to prove that the sequences $u_j = 4s_0 \delta^j$, 
$v_j = 4r_0 \delta^j$, and $w_j = 4h_0 \delta^j$, verify 
(\ref{c.3}) choosing  $r_0=cr$ and  $c=c(n,\rho,\tau, \varsigma)\ll 1$. 
We have $4s_0L_1 \le 1$ in view of (\ref{c.16a}). Moreover, 
$h_0 < \kappa \sigma_0^{\tau +1} < \kappa s_0 \le L_2^{-1-\varsigma}
s_0$, and we obtain $w_j L_2 \le u_j L_1$, and $w_j<
L_2^{-1-\varsigma}$. Finally, 
$r_0 = cr \le cL_2^{-1-\varsigma}< L_2^{-1-\varsigma}$, and 
$r_0L_2 < cL_1^{-\varsigma} \le s_0L_1$, choosing appropriately 
$c=c(n,\rho,\tau, \varsigma)\ll 1$.

It remains to show that
\begin{equation}
| P_0 |_{{\cal U}_{0}} \ \le \ \widetilde\varepsilon_0\ ,\quad 
|P_{j} - P_{j-1}|_{{\cal U}_{j}} \ \le\   
\widetilde \varepsilon_j,\ j\ge 1.  
                                     \label{c.18a}
\end{equation} 
for  $a\ll 1$. In view of (\ref{c.15}) we have 
$$
|P_0|_{{\cal U}_{0}}\le C_0\|P\| L_1^n = \widehat \varepsilon \kappa r_0 
C_0\frac{r}{r_0}L_1^{-N +n+2} a.
$$
On the other hand, using (\ref{c.16a}) we get  
$$
\begin{array}{rcl}
\sigma_0^{\tau +1} E_0 \, &=&\,     
c_1^{-1}  \sigma_0^{\tau +1} 
\exp \left( - A_0\,  {(L_1\sigma_{0})}^{ - {1 \over {\rho -1}}}
\right) \\ [0.3cm]
&=&  C'(n,\rho,\tau) 
L_1^{-\tau -1} (\ln (L_1 + e))^{-(\rho-1)(\tau + 1)}(L_1 + e)^{-M} 
\ge   C(n,\rho,\tau) L_1^{M-\tau -2}\, ,
\end{array}
$$
where $M = A_0(\rho)\sigma(n,\rho,\tau)^{-1/(\rho-1)}$. Now we fix 
$N = M+\tau
+n +4$ and chose 
$$ a = C(n,\rho,\tau)   C_0^{-1} \frac{r_0}{r} = 
 C(n,\rho,\tau)c(n,\rho,\tau, \varsigma) C_0^{-1}  .$$ 
Then $|P_0|_{{\cal U}_{0}}\le \widetilde\varepsilon_0$. 
Recall that $C_0$ comes from the Approximation lemma and the 
extension in the beginning of Sect. 3, hence,  $a$ is independent of 
$\kappa, L_1,L_2$, and $r$. 
Using (\ref{c.15}), we obtain   for each
$j\ge 1$ 
$$
\displaystyle{ |P_{j} - P_{j-1}|_{{\cal U}_{j}}  \le  
\widetilde\varepsilon_0 
\exp \left( - B_0\,  {\sigma_{j}}^{ - {1 \over {\rho -1}}}\right) \le 
\widetilde \varepsilon_j . }
$$

We are ready to make the iterations. We consider the real-analytic in 
${\cal U}_j$ Hamiltonian  
$H_j(\varphi,I;\omega) = N_0(I;\omega) + P_j(\varphi,I;\omega)$, where 
$ N_0(I;\omega) := \langle \omega,I\rangle$.  For any $j\in {\bf N}$,
we denote by ${\cal D}_j$ the class of  real-analytic 
diffeomorphisms 
${\cal F}_{j}: D_{j+1}\times  O_{j+1} \to D_{j}\times  O_{j}$
of the form  
\begin{equation}
{\cal F}_j(\theta,I;\omega)\ =\  (\Phi_j(\theta,I;\omega),\phi_j(\omega))\, ,\ 
\Phi_j(\theta,I;\omega)\ =\ (U_j(\theta;\omega),V_j(\theta,I;\omega))\, ,   
                                \label{c.19}
\end{equation}
where $\Phi_j(\theta,I;\omega)$ is affine  linear  with respect to $I$,  
and  the transformation  
$\Phi_j(.\, ,\, .\, ;\, \omega)$ 
is canonical  for  any fixed $\omega$. To simplify the notations we
denote the sup-norm in $D_{j}\times  O_{j}$ by $|\cdot |_j$ 
 instead of $|\cdot |_{s_j,r_j,h_j}$. Obviously, 
$D_{j}\times  O_{j} \subset {\cal U}_j$.

\begin{prop} \label{IterativeLemma}  
Suppose  $P_j$, $j\ge 0$,  is real-analytic on ${\cal U}_j$ with 
$$ 
| P_0 |_{{\cal U}_0}\ \le \ \widetilde\varepsilon_0\ ,\quad 
|P_{j} - P_{j-1}|_{{\cal U}_j} \ \le\   
\widetilde \varepsilon_j,\ j\ge 1.  
$$ 
Then for each $j\geq 0$ there exists a real analytic normal form 
$N_j(I;\omega) 
= e_j(\omega) + \langle \omega,I\rangle$ and a real analytic 
transformation  ${\cal F}^j$, where ${\cal F}^0 = {\rm
Id}$ and  
$$ 
{\cal F}^{j+1} = {\cal F}_0\circ\cdots\circ {\cal F}_{j}:\,  
D_{j+1}\times O_{j+1} \longrightarrow (D_0\times O_0) \cap {\cal
U}_{j}, \ j\ge 0, 
$$ 
with 
${\cal F}_j \in {\cal D}_j$  such that 
$H_j\circ {\cal F}^{j+1} = N_{j+1} + R_{j+1}$ and  
$\displaystyle | R_{j+1} |_{j+1}  \le \varepsilon_{j+1}$.  
Moreover, 
\begin{equation} 
|\overline{W}_j ({\cal F}_{j} - {\rm id} )|_{j+1} \ ,\ 
|\overline{W}_j (D {\cal F}_{j} - {\rm Id} )
\overline{W}_j^{-1}|_{j+1} \ <   
 \frac{C \varepsilon_j}{r_jh_j} , 
                         \label{c.20} 
\end{equation}
\begin{equation} 
|\overline{W}_0({\cal F}^{j+1} - {\cal F}^j)|_{j+1} \, 
< \frac{c \varepsilon_j}{r_jh_j} ,
                          \label{c.21} 
\end{equation} 
where $C =C(n,\rho)>0$ is the constant in Remark 3.3,
$c=c(n,\rho)>0$,   $D {\cal F}^{j}$
stands for the Jacobian of ${\cal F}^{j}$ with respect to
$(\theta,I,\omega)$, and 
$\overline{W}_j = {\rm diag}\, \left({\sigma_j}^{-1}{\rm Id}, 
{r_j}^{-1}{\rm Id},{h_j}^{-1}{\rm Id}\right)$.  
\end{prop} 
{\em Proof}. The proof is similar to that of the Iterative Lemma 
\cite{kn:Poe1}. First, applying the KAM step we find 
${\cal F}^{1} = {\cal F}_{0}$ such that 
$H_0 \circ {\cal F}_0 =  N_1 + R_{1}$, where $R_{1}$ is real analytic
in $D_1\times O_1$ and $|R_{1}|_{1} \leq \varepsilon_1$. 
By recurrence we define for any $j\ge 1$ the transformation 
${\cal F}^{j+1} = {\cal F}^{j}\circ {\cal
F}_{j}$, where ${\cal F}_j$  belongs to ${\cal D}_j$. By the inductive
assumption we have $H_{j-1} \circ {\cal F}^{j} = N_j + R_{j}$, where 
$N_j(I;\omega) = e_j(\omega) + \langle \omega, I\rangle$ is a
real-analytic normal form, $R_{j}$ is real analytic   in $D_j\times
O_j$, and  $ | R_{j} |_{j}  \leq  \varepsilon_{j}$. 
Then we write  
$$
H_j \circ {\cal F}^{j+1} = (N_0 + P_{j-1})\circ {\cal F}^{j+1} + 
(P_{j} - P_{j-1})\circ {\cal F}^{j+1} 
$$
$$
= \left[H_{j-1} \circ {\cal F}^{j}\right]\circ {\cal F}_{j} 
+ (P_{j} - P_{j-1})\circ {\cal F}^{j+1} 
$$
$$
= (N_j + R_{j})\circ {\cal F}_{j} + (P_{j} - P_{j-1})\circ {\cal
F}^{j+1} .
$$
We apply  Proposition \ref{kamlemma} to the
Hamiltonian $N_j + R_j$ which is
real-analytic in $D_j\times O_j$. In this way, using (\ref{c.14a}), 
we find a real-analytic map 
${\cal F}_{j}: 
D_{j+1} \times O_ {j+1} \to  D_j \times O_j$
which belongs to the class ${\cal D}_{j}$ and such that 
$(N_j + R_{j})\circ {\cal F}_{j} = N_{j+1} + R_{j+1,1}$, where 
$$ 
|R_{j+1,1}|_{j+1} \leq {1\over 2} 
\widehat\varepsilon  \kappa r_{j+1}\sigma_{j+1}^{\tau +1}
c_1^{1/2}E_j^{3/2} = {\varepsilon_{j+1} \over 2} .  
$$ 
Moreover,  ${\cal F}_{j}$ satisfies
(\ref{c.20}) in view of Remark \ref{estimates} and as in \cite{kn:Poe1}
we obtain  (\ref{c.21}). 
We are going to show that 
\begin{equation}
{\cal F}^{j+1} : D_{j+1} \times O_{j+1}\  \longrightarrow \ {\cal
U}_{j}\, . 
                                  \label{c.22}
\end{equation}
This inequality combined with  (\ref{c.18}) and (\ref{c.18a}) implies 
$$
|(P_{j} - P_{j-1})\circ {\cal F}^{j+1}|_{j+1} 
\leq |P_{j} - P_{j-1}|_{{\cal U}_{j}} \leq \widetilde \varepsilon_{j}
\leq {\varepsilon_{j+1} \over 2} ,   
$$
and  we obtain  
$H_j\circ {\cal F}^{j+1} = N_{j+1} + R_{j+1}$, where
$|R_{j+1}|_{j+1} \le \varepsilon_{j+1}$. 

To prove (\ref{c.22}) we note that  
$$
|\overline W_k \overline W_{k+1}^{-1}| = 
\sup \left\{s_{k+1}/s_k,r_{k+1}/r_k ,h_{k+1}/h_k
\right\} =  s_{k+1}/s_k = \delta, 
$$
since $r_{k+1}/r_k = E_k^{1/2} \ll \delta$, and $h_{k+1}/h_k \sim
\delta ^{(\tau +1)\rho (\rho-1)^{-1}} < \delta$ for any $k\in {\bf N}$. 
Then using (\ref{c.17a}) and (\ref{c.20}), 
we  estimate the Jacobian of ${\cal
F}^{j+1}$ in $D_{j+1}\times O_{j+1}$ as follows (see \cite{kn:Poe1})
$$
\begin{array}{rcl}
\left|\overline W_0 D{\cal F}^{j+1}\overline W_{j}^{-1}\right|_{j+1}
\,  &=& \, 
\left|\overline W_0 D ({\cal F}_{0}\circ \cdots \circ 
{\cal F}_{j})\overline W_{j}^{-1}\right|_{j+1}       \\[0.3cm]
&\le& \prod_{k=0}^{j-1}
\left(\left|\overline W_k D {\cal F}_{k}\overline W_{k}^{-1}\right|_{k+1}\, 
\left|\overline W_k \overline W_{k+1}^{-1}\right|\right)
\left|\overline W_j D {\cal F}_{j}\overline W_{j}^{-1}\right|_{j+1}\, 
 \\[0.3cm]
&\le& \delta^{j} \prod_{j=0}^\infty\, 
\left(1 +  \frac{C \varepsilon_j}{r_jh_j} \right)\ \le \  2\,   \delta^{j} \,
.
\end{array}
$$
Set $z=(\theta,I,\omega) =x +iy\in D_{j+1}\times  O_{j+1}$, where  $x$
and $y$ are the real and the imaginary part of $z$. Then 
$$
\begin{array}{rcl}
{\cal F}^{j+1}(x+iy) \, =\,  {\cal F}^{j+1}(x) +  
\overline W_0^{-1} T_{j+1}(x,y)\overline W_{j}\,  y \, ,\\ [0.3cm]
T_{j+1}(x,y) = \displaystyle { i\, \int_0^1 \,  \overline W_0 
D{\cal F}^{j+1}(x + ity) \overline W_{j}^{-1} \, dt \,  } .
\end{array} 
$$
Moreover,   
$|T_{j+1}(x,y)| \le 2\, \delta^{j}$ and using that  
$|\overline W_{j} \,  y| \le |\overline W_{j+1} \,  y| \le \sqrt 3$ we
get 
$$
|T_{j+1}(x,y)\overline W_{j}\,  y| \le 4 \, \delta^{j}\,  
,\quad x+iy\in D_{j+1}\times  O_{j+1} \, . 
$$
This implies 
${\cal F}^{j+1}(x+iy) \in {\cal U}_{j}$, since ${\cal F}^{j+1}(x)$ is
real, and we complete the proof of Proposition \ref{IterativeLemma}. 
\finishproof

We are going to prove suitable Gevrey estimates for ${\cal F}^{j}$. 
We set 
$$
\widetilde D_j =  \{(\theta,I)\in D_j :\, |{\rm Im}\, 
\theta|<s_j/2\}\ ,\quad 
\widetilde O_j = \{\omega\in {\bf C}^n:\ |\omega - \Omega_\kappa
|<h_j/2\} \, ,
$$ 
and we denote   ${\cal S}^j = {\cal F}^{j+1} - {\cal F}^j$. For any
 multi-indices $\alpha$ and $\beta$ and $m\in {\bf N}$ with
$|\beta|\le m$, we denote 
$$
R^m_{\omega '}\left(\partial_\theta^\alpha  \partial_\omega^\beta  
{\cal S}^j\right)(\theta, I, \omega) := 
  \partial_\theta^\alpha\partial^\beta_\omega {\cal S}^j(\theta, I, \omega)
-  \sum_{|\beta +\gamma|\le m}(\omega-\omega') ^{\gamma} 
\partial_\theta^\alpha\partial_\omega ^{\beta+\gamma}
{\cal S}^j(\theta, I, \omega' ) /\gamma ! .
$$
Recall that $\rho' = \rho(\tau +1) + 1$. 
\begin{lemma} \label{GevreyEstimates} 
Under the assumptions of Proposition \ref{IterativeLemma} we have
$$
\begin{array}{rcl}
|\overline{W}_0 \partial_\theta^\alpha \partial_\omega^\beta
{\cal S}^{j}(\theta,0,\omega)| \, &\le& \,  \widehat \varepsilon\,  A\, 
C^{|\alpha +\beta|}\, L_1^{|\alpha| +|\beta|(\tau + 1) + 1}
\kappa ^{-|\beta|}\, 
 \alpha !^{\, \rho} \beta  !^{\, \rho'}
E_j^{1/2} , \\ [0.3cm] 
&&(\theta,0,\omega) \in \widetilde D_{j+1}\times 
\widetilde O_{j+1}\, ,
\\ [0.5cm]
\displaystyle
|\overline{W}_0 
R^m_{\omega'} (\partial_\theta^\alpha\partial_\omega^\beta {\cal S}^{j}
(\theta,0,\omega))|\,
&\le& \, \widehat \varepsilon \,  A\,
C^{m+|\alpha|+1}L_1^{|\alpha| +(m+1)(\tau + 1) + 1}  
\kappa ^{-m-1}\, \\ [0.3cm]
&\times& \frac{|\omega-\omega'|^{m-|\beta|+1}}{(m-|\beta|+1)!}\, 
\alpha!\, ^\rho\,  (m+1) !\, ^{\rho'} E_j^{1/2}\, , \quad 
\theta\in {\bf T}^n,\ \omega, \omega '\in \Omega_\kappa\, ,
    \end{array} 
$$
for any 
$m\in {\bf N},\ \alpha, \beta  \in {\bf N}^n,\ |\beta| \le m$, 
where  the constants 
$A, C$ depend
only on $\tau$, $\rho$, $n$ and $\varsigma$.     
\end{lemma} 
{\em Proof}. Using (\ref{c.21}) and the Cauchy  estimate, 
we  evaluate $\partial_\theta^\alpha \partial_\omega^\beta
{\cal S}^{j}$  for any $j\ge 0$ and 
$|\alpha+ \beta|\ge 1$ in  
$\widetilde D_{j+1}\times  \widetilde O_{j+1}$. We have   
$$ 
M_{j,\alpha,\beta}\,  := 
\left| W_0 \partial_\theta^\alpha \partial_\omega^\beta {\cal S}^j\right| 
\, \le\,  c 
\frac{2^{|\alpha+\beta|}\, \alpha !\, \beta !\,
\varepsilon_j}{r_j h_j s_{j+1}^{|\alpha|} h_{j+1}^{|\beta|}} 
\,  =\,  c  \kappa \widehat \varepsilon 
\frac{2^{|\alpha+\beta|}\, \alpha !\,   \beta !\,
E_j\sigma_j^{\tau}}{{h_j s_{j+1}^{|\alpha|}} h_{j+1}^{|\beta|}}. 
$$ 
Recall that $s_j = 5(1-\delta)^{-1}\sigma_j = 
5(1-\delta)^{-1}A_0^{\rho -1} L_1^{-1} 
(B\sigma_j^{-\frac{1}{\rho-1}})^{-(\rho-1)}  $, and 
$$
h_j = {\kappa \over 2} K_j^{-\tau -1} = 
{\kappa \over 2} \sigma_j^{\tau +1}  x_j^{-\tau -1}\,  . 
$$
Then by  (\ref{c.17}) we get
$$
h_{j+1}^{-1} \le   \kappa^{-1} \widetilde C_0  L_1^{\tau +1} 
\left(B  \sigma_{j}^{-\frac{1}{\rho -1}}\right)^{\rho(\tau +1)}.
$$
where  $\widetilde C_0$ depends only on $\tau$ and $\rho$. 
This implies 
\[
\begin{array}{rcl}
M_{j,\alpha,\beta} \, &\le&   \,\widehat \varepsilon\, A_1
 C_1^{|\alpha + \beta|} L_1^{|\alpha| +|\beta|(\tau + 1) + 1}
\kappa^{-|\beta|} 
\alpha !\, \beta  !\, \\ [0.5cm]
&\times& \, \left(B\sigma_{j}^{-\frac{1}{\rho-1}}\right)^
{(\rho-1)(|\alpha| -\tau) + \rho (\tau +1)(|\beta|+1)}
\, \exp \left(-B \sigma_{j}^{-\frac{1}{\rho-1}}\right) \, , 
\end{array}
\]
where $A_1, C_1$  depend 
only on $\tau$ and $\rho$. 
Then  we obtain  
\begin{equation} 
M_{j,\alpha,\beta} \, \le  \, \widehat \varepsilon A \,
 C^{|\alpha + \beta|} L_1^{|\alpha| +|\beta|(\tau + 1) + 1}
\kappa^{-|\beta|} 
\alpha !^{\, \rho}\, \beta  ! ^{\, \rho (\tau +1) +1}  \, E_j^{1/2} 
\, . 
                                 \label{c.24}
\end{equation} 
where   $A, C$ depend 
only on $\tau$ and $\rho$. 

We are going to prove the second estimate for $\omega, \omega' \in
\Omega_\kappa$.  First we suppose 
that $|\omega' - \omega| \le h_{j+1}/8$. 
Expanding  the analytic in $O_{j+1}$ function    
$\omega \to \partial^\alpha_\omega {\cal S}^{j}(\theta,0, \omega)$,
$\theta \in {\bf T}^n$,  
in Taylor series with respect to $\omega$ at $\omega'$, 
and using as above the Cauchy estimate for  
$M_{j,\alpha,\gamma}$, we  evaluate    
$$
L_{j,\alpha,\beta}^m\,  := \, 
|W_0(R_{\omega'}^m\, \partial_\theta^\alpha   
\partial_\omega^\beta {\cal S}^j)(\theta,0;\omega))|\ ,\quad 
\theta\in  {\bf T}^n,\ \omega,\, \omega' \in \Omega_\kappa . 
$$
For $|\beta|\le  m + 1$ we have 
$$
\frac{(\beta + \gamma)!}{\gamma !} \le 2 ^ {|\beta + \gamma|} \beta !
\le 2 ^ {|\beta + \gamma|} \frac{(m+1)!}{(m-|\beta|+1)!} .
$$
Then we obtain   as above 
\[
\begin{array}{rcl}
L_{j,\alpha,\beta}^m\ &\le &  \ 
\displaystyle{ \sum_{|\gamma|\ge  m -|\beta| +1}\, 
|\omega' - \omega|^{|\gamma|}\, 
M_{j,\alpha ,\beta+\gamma}(\theta,0,\omega')/\gamma !}
 \\ [0.5cm]
 &\le& \displaystyle{ c \, \alpha ! (m+1)! \, \frac{|\omega' -
\omega|^{m-|\beta|+1}}{(m-|\beta|+1)!}\,  
\frac{4^{|\alpha|+m+1} \varepsilon_j}
{r_j h_j s_{j+1}^{|\alpha|} h_{j+1}^{m+1}} \,
\sum_{|\beta + \gamma|\ge  m  +1}\, 
\left(4 |\omega' - \omega|h_{j+1}^{-1}\right)^{|\beta +\gamma |-m-1} },
\end{array}
\]
and we get 
\[
\begin{array}{rcl}
L_{j,\alpha,\beta}^m\ &\le &  \displaystyle {
2  c \,\alpha ! (m+1)!  \, 
\frac{|\omega' - \omega|^{m-|\beta|+1}}
{(m-|\beta|+1)!}\, 
\frac{4^{|\alpha|+m+1} \varepsilon_j}
{r_j h_j s_{j+1}^{|\alpha|} h_{j+1}^{m+1}} \,}
\\ [0.5cm]
 &\le& \displaystyle{  \widehat \varepsilon\,  A\,
C ^{|\alpha| + m+1} L_1^{|\alpha| +(m+1)(\tau + 1) + 1}  \kappa^{-m-1}\,  
\frac{|\omega' - \omega|^{m -|\beta|+1}}{(m-|\beta|+1)!}\, 
\alpha !^{\rho}  (m+1)!\, ^{ \rho (\tau +1) +1}  E_j^{1/2}. }
\end{array}
\]
where  $A, C$ depend
only on $\tau$, $\rho$ and $n$.  For  
$|\omega' - \omega| \ge h_{j+1}/8$ we obtain the same inequality, 
estimating $L_{j,\alpha,\beta}^m$ term by term and using (\ref{c.24}). 
This proves the lemma. \finishproof

According to Proposition \ref{IterativeLemma} and 
Lemma \ref{GevreyEstimates},  the
limit 
$$
\partial_\theta^\alpha {\cal H}^{\beta}(\theta,\omega)  := 
\lim_{j\to \infty} \partial_\theta^\alpha
\partial_\omega^\beta \left[ {\cal F}^j(\theta,0,\omega) -
(\theta,0,\omega) \right]\, ,\quad 
(\theta,\omega)\in {\bf T}^n\times \Omega_\kappa\, ,
$$
exists for each $\alpha,\beta  \in {\bf N}$ and it is uniform since 
$\sum E_j^{1/2} < \infty$.  
Moreover, the partial derivatives of $ \partial_\theta ^\alpha ({\cal
H}^{\beta}) =  \partial_\theta ^\alpha {\cal
H}^{\beta}$ exist and they are continuous on ${\bf T}^n\times
\Omega_\kappa$. 
Consider the jet
${\cal H}=\left(\partial_\theta^\alpha{\cal H}^{\beta}\right)$, $\alpha,\beta
\in {\bf N}^n$, 
 of  continuous functions
$\partial_\theta^\alpha{\cal H}^{\beta} : 
{\bf T}^n\times \Omega_\kappa \to {\bf T}^n\times D \times  \Omega$,
and  set 
$$
\left(R^m_{\omega '} 
\partial_\theta^\alpha {\cal H}\right)_\beta(\theta,\omega) := 
  \partial_\theta^\alpha{\cal H}^\beta(\theta,\omega)
-  \sum_{|\beta +\gamma|\le m}(\omega-\omega') ^{\gamma} 
\partial_\theta^\alpha
{\cal H}^{\beta+\gamma}(\theta, \omega' ) /\gamma ! .
$$
In view of Lemma \ref{GevreyEstimates}, we have
\begin{equation}  
\begin{array}{rcll} 
|\overline{W}_0 \partial_\theta^\alpha 
{\cal H}^{\beta}(\theta,\omega)|  &\le&  \widehat
\varepsilon \, A\, L_1\, 
 (CL_1)^{|\alpha|} \, (CL_1^{\tau + 1}\kappa^{-1})^{|\beta|}   
\alpha !^{\, \rho} \beta  !^{\, \rho'}
\nonumber \\ [0.3cm]
|\overline{W}_0\left(R^m_{\omega '} 
\partial_\theta^\alpha {\cal H}\right)_\beta(\theta,\omega)|
&\le&  \widehat \varepsilon\,  A\,  L_1\,  
(CL_1)^{|\alpha|} \,  (CL_1^{\tau + 1}\kappa^{-1})^{m+1} 
\displaystyle \frac{|\omega-\omega'|^{m-|\beta|+1}}{(m-|\beta|+1)!}\, 
\alpha!\, ^\rho\,  (m+1) !^{\, \rho'} \, 
                                         \label{c.25} 
\end{array}
\end{equation} 
for each  $\alpha$,  and $\beta$ satisfying 
$0\le |\beta|\le m$, and 
  $\theta \in  {\bf T}^n$,  
$\omega,\omega'\in \Omega_\kappa$, where  $A$ and $C$ depend only on 
$\tau$, $\rho$,  $\rho'$, and $n$. 
We are going to extend ${\cal H}$ to a Gevrey function on 
${\bf T}^n\times \Omega$.

\subsection{\it Whitney extension in Gevrey classes}
\label{Sec:Whitney}
Let $K$ be a compact in ${\bf R}^n$ and $\rho \ge 1$, $ \rho' >1$. 
We consider a jet $(f^{\beta})$, $\beta\in {\bf N}^n$, of functions 
$f^{\beta}: {\bf T}^n\times K \to {\bf R}$, such that for each 
$\alpha\in {\bf N}^n$ the partial derivative 
$\partial_\theta^\alpha f^{\beta}$ exists,  it 
is  continuous  on ${\bf T}^n\times K$, and there are positive 
constants $A$, $C_1$ and  $C_2$ such that 
\begin{equation}  
\begin{array}{rcll} 
|\partial_\theta^\alpha f^{\beta}(\theta,\omega)| \, &\le& \, \, A\, 
C_1^{|\alpha|} C_2^{|\beta|}\, \alpha !^{\, \rho} \beta  !^{\, \rho'},
\nonumber \\ [0.3cm]
|\left(R^m_{\omega '} 
\partial_\theta^\alpha f\right)_\beta(\theta,\omega)|\,
&\le& \, A\, 
C_1^{|\alpha|} C_2^{m+1}\,
\displaystyle\frac{|\omega-\omega'|^{m-|\beta|+1}}{(m-|\beta|+1)!}\, 
\alpha!\, ^\rho\,  (m+1) !^{\, \rho'} \, . 
\end{array}
                              \label{c.26} 
\end{equation}
\begin{Theorem} \label{Whitney}
There exist positive constants $A_0$ and $ C_0$   
and for any compact set $K$ and
any  jet $f=(f^\beta)$, $\beta\in {\bf N}^n$, satisfying (\ref{c.26}) there
exists  $\widetilde f\in {\cal G}^{\rho,\rho'}({\bf T}^n\times
{\bf R}^n)$ such that $\partial_\theta^\alpha\partial_\omega ^{\beta}
\widetilde f
= \partial_\theta^\alpha f^\beta$ on $K$ for each $\alpha$, $\beta$,
and 
$$
|\partial_\theta^\alpha\partial_\omega ^{\beta} \widetilde 
f(\theta,\omega)| \, 
\le \, A\,  A_0 \,  \max (C_1,1)\,  
(C_0 C_1)^{|\alpha|+1}\,(C_0 C_2)^{|\beta|}\,  
\alpha !^{\, \rho} \beta  !^{\, \rho'}\, .
$$
\end{Theorem}
{\em Remark}.  We point out that the  positive 
constants $A_0$ and $C_0$ do not depend on the
the jet $f$, on compact set $K$ nor on the constants $A,\, C_1$ and $C_2$.

\noindent {\em Proof.} \hspace{2mm} 
Consider the Fourier coefficients 
$$ 
f_k^\beta(\omega)\ =\ (2\pi)^{-n}\, \int_{{\bf T}^n}\, 
e^{-i\langle k,\varphi\rangle}\, 
f^\beta (\varphi,\omega)d\varphi\ ,\ k\in {\bf Z}^n\ , 
$$ 
and denote by ${\cal A}$ the set of all Whitney jets $g_k =
\left(g_k^\beta\right)$, $\beta\in {\bf N}^n$, where 
$g_k^\beta = e^{r|k|^{1/\rho}}\, f_k^\beta(\omega)$ and $k\in {\bf
Z}^n$, $|k| = |k_1| + \cdots + |k_n|$. Choosing   
  $r= c_0 C_1^{-1/\rho}$ with $0<c_0 = c_0(n,\rho) \ll 1$, 
we are going to show that $(g_k^\beta)\in {\cal A}$
satisfy 
\begin{equation}  
\begin{array}{rcll} 
\left|g_k^\beta (\omega)\right|\, &\le& \,  
A_2 C_2 ^{|\beta|}\,  \beta !^{\, \rho}\ ,\ \omega\in K\, 
\\ [0.3cm]
|\left(R^m_{\omega '} 
g_k \right)_\beta(\omega)|\, &\le& \, A_2
C_2^{m+1}\, 
\displaystyle\frac{|\omega-\omega'|^{m-|\beta|+1}}{(m-|\beta|+1)!}\, 
  (m+1) !^{\, \rho'} \, , 
\end{array}
                               \label{c.27}
\end{equation}
where $A_2 = 2 A  \max (C_1,1)\,  $. 
We decompose $j=\rho \ell_j + q_j$, $j,\ell\in {\bf N}$, $0\le
q_j<\rho$. For any  
$k=(k_1,\ldots,k_n)\in {\bf Z}^n$ we have 
$|k| \le  n \max_{1\le i \le n} |k_i|= n|k_p|$ for some 
$p$. Then integrating  by parts and 
using (\ref{c.26}), we estimate 
$$
\frac{r^j|k|^{j/\rho}}{j!} |f_k^\beta (\omega)| 
\le  A \frac{r^j}{j!} n^{\ell_j+1} C_1^{\ell_j+1} C_2^{|\beta|} 
(\ell_j+1)!^\rho \beta ! \, ^{\rho'}
\le  2^{-j}  A  \max (C_1,1)\,   C_2^{|\beta|}
\beta ! \, ^{\rho'} ,
$$
choosing $r= c_0 C_1^{-1/\rho}$ with $0<c_0 = c_0(n,\rho) \ll 1$. This
proves the first part of (\ref{c.27}). To prove the second part, we
notice that $\left(R^m_{\omega '} 
f_k \right)_\beta(\omega)$ is just the Fourier coefficient of 
$\left(R^m_{\omega '} 
\partial_\theta^\alpha f\right)_\beta(\theta,\omega)$ corresponding to
$k$. 
Now we use a  variant of the Whitney extension theorem due to Bruna 
(see Theorem 3.1, \cite{kn:Bru}). 
\begin{Theorem} \label{Bruna}
For any compact set $K$ and
a jet $g=(g^\beta)$, $\beta\in {\bf N}^n$, satisfying (\ref{c.27}) on
$K$,  there
is $\widetilde g\in {\cal G}^{\rho'}({\bf R}^n)$ such that
$\partial_\omega ^{\beta} \widetilde g = g^\beta$ on $K$ for any $\beta$, and 
$$
|\partial_\omega ^{\beta} \widetilde  g(\omega)| \, 
\le \, \, A_0 A_2\, 
(C_0 C_2)^{|\beta|}\, \beta  !^{\, \rho'}\, .
$$
Moreover, the positive constants $A_0$ and $C_0$ do not depend on the
jet $g$,
on the
compact set $K$ nor on the constants $A_2,\, C_2$. 
\end{Theorem}
The proof of Theorem \ref{Bruna} is given in \cite{kn:Bru}. Here we
only indicate that the constants $A_0,\, C_0 >0$ do not depend on  
the compact set $K$ nor on $A_2$ and $C_2$. This  follows from the proof
of Theorem 3.1,  \cite{kn:Bru}. 
More precisely, setting  $f^\beta(\omega)=A_2^{-1}C_2^{-|\beta|}g^\beta
(C_2^{-1}\omega)$, $C_2^{-1}\omega \in K$, we obtain 
$|f^\beta(\omega)| \le \beta  !^{\, \rho'}$ and
$$
|\left(R^m_{\omega '} f \right)_\beta(\omega)|\, 
= \, A_2^{-1}C_2^{-|\beta|} 
\left|\left(R^m_{C_2^{-1}\omega '} g \right)_\beta(C_2^{-1}\omega)\right|\, 
\le \, 
\frac{|\omega-\omega'|^{m-|\beta|+1}}{(m-|\beta|+1)!}\, 
  (m+1) !^{\, \rho'} \,   
$$
for $C_2^{-1}\omega, C_2^{-1}\omega' \in K$. 
Hence, we can suppose that (12) and (13), \cite{kn:Bru}, hold with
$\varepsilon =A=1$ ($K$ is  scaled to
another compact still denoted by $K$). 
Then it is easy to see that 
the constants $A$ and $\varepsilon$ in (19),
\cite{kn:Bru}, are independent of $g$ and $K$. Moreover, the different
constants in Lemma 3.2 and 3.3, \cite{kn:Bru}, do not depend on $g$
and $K$,
and we obtain that the constants in (26), \cite{kn:Bru}, are
independent of $g$ and $K$. Scaling back by $C_2$  we obtain the desired
estimates with constants $A_0$ and $C_0$ independent of $g$, $K$, $A_2$ and
$C_2$.  \finishproof

Applying Theorem \ref{Bruna} to the family of jets ${\cal A}$, we
obtain a family of functions $\widetilde g_{k}\in {\cal G}^{\rho '}$ 
such that 
$\partial_\omega ^\beta \widetilde g_{k} = g_{k}^\beta$ on $K$ for each
$\beta$, and 
$$ 
\left|\partial_\omega^\beta\, \widetilde g_{k}(\omega)\right|\ 
\leq 
A_0 A \max(C_1,1) (C_0 C_2)^{|\beta|}\,  
\beta !^{\, \rho'}\ ,\ \forall\, \omega \in
{\bf R}^n\, , \ 
k\in {\bf Z}^n\, ,\ \beta\in {\bf N}^n\, . 
$$ 
Now it is easy to see that  the function 
$$
\widetilde f (\theta,\omega)\ =\ \sum_{k\in {\bf Z}^n}\:  
e^{i\langle k,\theta\rangle-r|k|^{1/\rho}}\, 
\widetilde g_{k}(\omega)  
$$
satisfies the  requirements of Theorem \ref{Bruna}. 
\hfill $\Box$ \vspace{5mm}

\subsection{\it Proof of Theorem 2.1.}
\label{Sec:Proof}
\vspace{0.5cm}
\noindent
Using (\ref{c.25}), we extend the jet ${\cal H}$ to a Gevrey function 
${\cal H}=({\cal H}_1,{\cal H}_2,{\cal H}_3): 
{\bf T}^n\times \Omega \to  {\bf T}^n\times D \times
\Omega$. We have 
$$
|\overline{W}_0 \partial_\theta^\alpha \partial_\omega^\beta
{\cal H}(\theta,\omega)| \, \le \, 
 \widehat\varepsilon \, A\, L_1^2\, 
 (CL_1)^{|\alpha|} \, (CL_1^{\tau + 1}\kappa^{-1})^{|\beta|}   
\alpha !^{\, \rho} \beta  !^{\, \rho(\tau' +1) +1}\, ,
$$
where $\widehat\varepsilon = \|P\|L_1^{N-2}(a\kappa r )^{-1}$ and   
the positive constants $A$ and $C$  are independent  of 
$L_1$, $L_2$, $\kappa$, $r$,  and $\Omega \subset
B_{\bar R}(0)$.  We  set ${\cal F} = (\Phi,\phi)$, $\Phi =(U,V)$, where 
$U(\theta,\omega) = {\cal H}_1(\theta,\omega) +\theta$, 
$V(\theta,\omega) = {\cal H}_2(\theta,\omega)$, and $\phi(\omega) = 
{\cal H}_3(\omega) +\omega$.  
Recall that  $r_0=cr$, where $c=c(n,\tau,\rho,\varsigma) >0$ is fixed  in
Sect. 3.3. 
On the other hand,  $h_0 \le \kappa
\sigma_0^{\tau +1} <  \kappa$, and   we obtain
\[
\begin{array}{lrc}
\left|\partial_\theta^{\alpha}\partial_\omega^{\beta}
(U(\theta;\omega) - \theta)\right|\,  
+ \, r^{-1} 
\left|\partial_\theta^{\alpha}\partial_\omega^{\beta}V(\theta;\omega)\right|  
+ \kappa^{-1} \left|\partial_\omega^{\beta}
(\phi(\omega) - \omega)\right|
\\    [0.3cm]  
\le\ \displaystyle  
A  \frac{\| P \|L_1^{N }}{\kappa r }    \, 
 (CL_1)^{|\alpha|} \, (CL_1^{\tau + 1}\kappa^{-1})^{|\beta|}  
\alpha !\, ^{\rho } \beta !\, ^{\rho '}        
\end{array}
\]
for some positive constants $A$ and $C$ as above. Choosing $\epsilon
<1/A$ in (\ref{b.3}) we obtain 
$|V(\theta,\omega)| \le A \| P \|L_1^{N }\kappa ^{-1}
\le A \epsilon r < r \le R$. In the same way we
get $\phi(\omega)\in \Omega$ for $\omega\in \Omega_\kappa$. 
This proves the
estimates  in Theorem 2.1. 
As in Sect. 5.d, \cite{kn:Poe1}, we obtain that
$\left| X_{H_j} \circ {\cal F}^j - D \Phi^j \cdot X_N\right| \le 
\frac{c\varepsilon_j}{r_jh_j}$
on ${\bf T}^n \times \{0\}\times \Omega_\kappa$ for all $j\ge 0$,
where $ X_{H_j}$ and $X_N$ stand for the Hamiltonian vector fields of 
$H_j(\theta,I;\omega)$ and $N=\langle \omega,I \rangle $,
respectively. On the other hand, $\nabla H_j$ converges uniformly to 
$\nabla H$ as $j\to \infty$ in view of Proposition 3.1, hence, 
$$
X_{H(\cdot;\phi(\omega))} \circ \Phi = D \Phi \cdot X_N
$$
on ${\bf T}^n \times \{0\}\times \Omega_\kappa$. Then 
$\{\Phi(\theta;\omega):\, \theta\in {\bf T}^n\}$ is an embedded invariant
torus 
of the Hamiltonian $H(\theta,I;\phi(\omega))$
with frequency $\omega \in \Omega_\kappa$. It is Lagrangian by
construction (see also \cite{kn:Her}, Sect. I.3.2). 
This completes the proof of Theorem 1.1.  \finishproof

\subsection{\it Real analytic hamiltonians. }
\label{Sec:Analytic}
\vspace{0.5cm}
\noindent
Consider a real analytic Hamiltonian 
$P \in {\cal G}^1_{ L_1,L_2,L_2}({\bf T}^n\times B \times
\Omega)$, $B=B_R(0)$,  with norm 
$$
\|P\|_{L_1,L_2} = \sup  
\left(|\partial_\theta^\alpha \partial_I^\beta \partial_\omega^\gamma 
P(\theta, I;\omega)| L_1^{-|\alpha|} L_2^{-|\beta|-|\gamma|}
(\alpha ! \beta ! \gamma !)^{-1}\right) < +\infty . 
$$
Then $P$ can be extended as an analytic Hamiltonian in 
$D_{s,r}\times O_h$, where $s= (2L_1)^{-1}$ and $r=h= (2L_2)^{-1}$,
and  with sup-norm satisfying $\|P\|_{2L_1,2L_2} \le \|P\|_{s,r,h}
\le C\|P\|_{L_1,L_2}$, where $C=C(n) >0$. 
Fix $\tau'>\tau>n-1$. We can slightly improve Theorem 3.1
\cite{kn:Pop}. Given $s>0$ we denote by  ${\cal U}_s$ the set of all
$\theta \in {\bf C}^n /2\pi {\bf Z}^n$ such that $|{\rm Im}\, \theta |
< s$. 
\begin{Theorem} \label{analyticKAMtheorem}
Suppose that $H$ is given by  (\ref{b.1}) where  $P \in {\cal
G}^1_{ L_1,L_2,L_2}({\bf T}^n\times 
B_R(0) \times  \Omega)$. Fix  
  $\kappa >0$ and $r>0$ such that 
$\kappa, r <  L_2^{-1}$ and $r \le R$.   Then there is $0<s_0\le s$
and 
$\epsilon>0$ both independent of $\kappa$,  $L_2$, $r$, $R$, and of 
$\Omega \subset \Omega_0$, 
such that if  
$\| P \| \leq\ \epsilon \kappa r $
then   
there  exist maps $\phi \in {\cal G}^{\tau' +2}(\Omega, \Omega)$ and 
$\Phi = (U,V) 
\in {\cal G}^{1,\tau ' +2}({\cal U}_{s_0/2}\times\Omega,\,  
{\cal U}_{s_0}\times B_R(0))$,
satisfying (i), Theorem 2.1,  with $\rho=1$. Moreover, 
there exist $A, C >0$, independent of 
$\kappa$,    $L_2$, $r$,  and $\Omega$, such that  
\[
\left|\partial_\omega^{\beta}
(U(\theta;\omega) - \theta)\right|\,  + \,
r ^{-1}
\left|\partial_\omega^{\beta} V(\theta;\omega)\right|  
+ \kappa^{-1} \left|\partial_\omega^{\beta}
(\phi(\omega) - \omega)\right|
\leq\  A\,  C^{|\beta|}\,  
\kappa^{-|\beta|} \, \beta !^{\, \tau' + 2}\,  \displaystyle 
 \frac{\| P \|} {\kappa r}       
\]
 uniformly in  $(\theta,\omega)\in {\cal U}_{s_0/2}\times \Omega$ and 
for any $\beta\in {\bf N}^n$.    
\end{Theorem} 
{\em Proof}. Fix $\rho = (\tau'-\tau)(\tau +1)^{-1} +1$ and set as 
above $\sigma_j = \sigma_0 \delta^j$,  and 
$\displaystyle s_{j+1} = s_j - 5\sigma_j$, $j\ge 0$, 
where $ s_0(1-\delta) = 20 \sigma_0$ (then $s_j\to 3s_0/4$) and $s_0
\le s = (2L_1)^{-1}$.  
As above we define 
$E_j := c_1^{-1}\exp ( - \sigma_{j}^{ - {1 \over {\rho -1}}})$
and set 
$\varepsilon_j =\widehat \varepsilon  \kappa r_j \sigma_j^{\tau +1}
E_j$, where $\widehat \varepsilon = \| P\|_{s,r,h} (a\kappa r)^{-1}$. 
Choose $r_0=r$ and $h_0=h$ and define $r_j$, $x_j$, $K_j$ and $h_j$ as
above. 
We consider $H=\langle \omega, I\rangle + P$ in $D_j\times O_j$. 
We fix $\sigma_0 =\sigma_0(n,\rho,\tau, s) \le s(1-\delta)/20$,
so that a), b), c) in Proposition 3.1  and (\ref{c.14}) hold for any
$j\ge 0$. 
Next we choose 
$a = a(n,\rho,\tau, L_1) \le  
c_1^{-1} \sigma_{0}^{ - \tau -1} \exp ( - \sigma_{0}^{ - {1 \over
{\rho -1}}})$. Then 
$|P|_0 \le \varepsilon _0$, and we can apply Proposition
3.5. Moreover, as in Lemma 3.6 we obtain with $\rho' =\rho(\tau +1) +1
= \tau' + 2$ 
$$
\begin{array}{rcl}
|\overline{W}_0 \partial_\omega^\beta
{\cal S}^{j}(\theta,0,\omega)| \, &\le& \,  \widehat \varepsilon\,  A\, 
C^{|\beta|}\, 
\kappa ^{-|\beta|}\, 
\beta  !^{\, \rho'}
E_j^{1/2} ,\  (\theta,0,\omega) \in {\cal U}_{s_0/2}\times 
\widetilde O_{j+1}\, ,
\\ [0.5cm]
\displaystyle
|\overline{W}_0 
R^m_{\omega'} (\partial_\omega^\beta {\cal S}^{j}
(\theta,0,\omega))|\,
&\le& \, \widehat \varepsilon \,  A\,
C^{m+1} 
\kappa ^{-m-1}\, \\ [0.3cm]
&\times& \frac{|\omega-\omega'|^{m-|\beta|+1}}{(m-|\beta|+1)!}\, 
(m+1) !\, ^{\rho'} E_j^{1/2}\, , \quad 
\theta\in {\cal U}_{s_0/2},\ \omega, \omega '\in \Omega_\kappa\, ,
    \end{array} 
$$
for any 
$m\in {\bf N},\ \beta  \in {\bf N}^n,\ |\beta| \le m$, 
where  the constants 
$A, C$ depend
only on $\tau$, $\rho$, $n$ and $L_1$. We complete the proof of the
theorem as above. \finishproof

\section{\it Proof of Theorem \ref{KAMtheorem} and Corollary 1.2. } 
\setcounter{equation}{0} 

{\em Proof of Theorem \ref{KAMtheorem}}. 
Set $r=R=\kappa \sqrt{\epsilon_H}$. Then (\ref{b.2}) implies 
$\|P\| \le (A + 1)\kappa r \sqrt{\epsilon_H}  $ 
and we can apply Theorem 2.1. 
Consider the map $\overline \Phi : {\bf T}^n \times \Omega
\to {\bf T}^n \times D$ given by 
$$
 \overline\Phi(\theta,\omega) = 
(U(\theta;\omega), \nabla g^0(\phi(\omega)) + V(\theta;\omega)),
$$
where $U$ and $V$ are obtained in  Theorem \ref{KAMtheorem2}. 
We have   $H(\theta,I; \phi(\omega)) = 
H(\theta,  \nabla g^0(\phi(\omega)) + I)$, $I\in B_R(0)$, in  
the notations of Sect. 2.1.  Then 
Theorem 2.1, (i),  implies that 
$\Lambda_\omega = \{\overline  \Phi(\theta,\omega):\, \theta\in
{\bf T}^n\}$, $\omega \in \Omega_\kappa$,
is an embedded Lagrangian invariant torus of the
Hamiltonian $H$  for each 
$\omega\in \Omega_\kappa$ with frequency $\omega$. The
corresponding  estimates for $\overline \Phi$ follow directly from
those in Theorem 2.1.  \finishproof

\noindent
{\em Proof of Corollary 1.2}. The proof is close to that of Theorem
2.1, \cite{kn:Pop} and we will be concerned  mainly with 
the corresponding Gevrey
estimates. Let $\varepsilon_H \le \epsilon L_1^{-N-2}$. 
Then in the notations of (ii), Theorem 1.1, we get 
$A C_1  L_1^{N/2} \sqrt{\varepsilon_H} \le AC\epsilon$. 
Choosing $\epsilon$ small enough and 
using Proposition A.2   as well as  (ii), Theorem 1.1, we obtain a
solution $\theta = \theta(\varphi, \omega)$ 
the equation $\overline U(\theta,\omega)=\varphi$ such that  
$ \theta(\varphi, \omega) -\varphi$ satisfies the same  Gevrey
estimates as  $\overline U$. 
Set $F(\varphi,\omega)= \overline V(\theta(\varphi, \omega),\omega)$. 
We have  $\Lambda_\omega = \{(\varphi,F(\varphi,\omega)):\, \varphi\in
{\bf T}^n\}$ for each $\omega\in \Omega_\kappa$. Moreover, 
by Proposition A.4  we obtain
$$
\left|\partial_\varphi^{\alpha}\partial_\omega^{\beta}
(F(\varphi,\omega) - \nabla g^0(\omega))\right|\,  
\le\  \kappa  A  C_1^{|\alpha|} (C_2\kappa^{-1})^{|\beta|}  
\alpha !\, ^{\rho } \beta !\, ^{\rho(\tau + 1) +1}\,   
L_1^{N/2} \sqrt{\varepsilon_H} 
$$
uniformly in  $(\theta,\omega)\in {\bf T}^n\times \Omega$. Hereafter, 
$A$, $ C_1=CL_1$, and $C_2=CL_1^{\tau +1}$ 
are positive constants  as in Theorem 1.1. 
Denote by $p:{\bf R}^n \to  {\bf T}^n$ the natural
projection. As in Lemma 2.2 , \cite{kn:Pop}, we shall find 
$\psi\in  {\cal G}^{\rho,\rho'}( {\bf R}^n\times \Omega)$ and 
$R\in  {\cal G}^{\rho'}(\Omega)$ such that 
$Q(x,\omega):= \psi(x,\omega)-\langle x,R(\omega)\rangle$ is $2\pi$
periodic with respect to $x$ and 
\begin{itemize}
\item[(i)] $\forall (x,\omega)\in{\bf R}^n\times \Omega_\kappa$, 
 $\nabla_x\psi(x,\omega)= F(p(x),\omega)$,
\item[(ii)]   
$\left| \partial_x^\alpha \partial_\omega^\beta 
Q(x,\omega) \right|\, 
+\, \left|\partial_I^\beta (R(\omega) - \nabla g^0(\omega))\right|\ 
 \le\  \kappa  A  C_1^{|\alpha|} (C_2\kappa^{-1})^{|\beta|}  
\alpha !\, ^{\rho } \beta !\, ^{\rho(\tau + 1) +1}\,
L_1^{N/2} \sqrt{\varepsilon_H} ,$
\end{itemize} 
for $(x,\omega)\in{\bf R}^n\times \Omega$. 
To obtain $\psi$ we consider the function
$$ 
\widetilde\psi (x,\omega)\  =\  \int_{\gamma_x}^{}\, \sigma\ 
=\ \int_{0}^{1}\, \langle F(p(tx),\omega),\, x\rangle\, dt, \quad 
 (x,\omega)\in {\bf R}^n\times \Omega\, ,
$$ 
where $\gamma_x = \{(tx, F(p(tx),\omega)): 0\leq t\leq 1\}$ and $\sigma =\xi dx$ is the canonical one-form on $T^\ast{\bf
R}^n$. Then 
$\widetilde\psi(x, \omega) - \langle \nabla g^0(\omega),x \rangle$
satisfies the Gevrey
estimates (ii) in $[0,4\pi]^n\times \Omega$.  We set 
$2\pi  R_j (\omega) = \widetilde\psi(2\pi e_j,\omega)$, 
$\omega\in \Omega$, 
$\{e_j\}$ being an unitary  basis in ${\bf R}^n$. Then $R -  \nabla
g^0$ satisfies (ii) in $\Omega$.  
Since $\Lambda_\omega$, $\omega\in \Omega_\kappa$,  is Lagrangian, 
we obtain as in \cite{kn:Pop} 
that for such $\omega$ the function 
$\nabla_x \widetilde \psi(x,\omega)$ is $2\pi$ periodic with
respect to $x$ and
$\widetilde\psi (x+2\pi m,\omega) - \widetilde\psi (x,\omega) = 
\langle 2\pi m, R(\omega)\rangle$ for $m\in {\bf Z}^n$. 
Consider the function 
$\widetilde Q(x,\omega) = 
\widetilde\psi(x,\omega) - \langle x, R(\omega) 
\rangle$. It satisfies the Gevrey estimates (ii) in 
$[0,4\pi]^n\times \Omega$, and it is $2\pi$ periodic with respect to
$x$ for $\omega\in \Omega_\kappa$. 
We are going to average  $\widetilde Q$ on ${\bf T}^n$. Let 
$f\in {\cal G}^\rho_C({\bf R}^n)$ with ${\rm supp}\, f \subset [\pi/2,
7\pi /2]^n$, where $C>0$ is a positive constant, and such that 
$\sum_{k\in {\bf Z}^n} f(x - 2\pi k) = 1$ for each $x\in {\bf R}^n$. 
Consider the function $Q(x,\omega) =
\sum_{k\in {\bf Z}^n} (f \widetilde Q)(x - 2\pi k,\omega)$.  
It is $2\pi$-periodic with respect to $x$ by
construction, it belongs to
${\cal G}^{\rho,\rho'}({\bf R}^n\times \Omega)$, and 
$Q(x,\omega) = \widetilde Q(x,\omega)$ for $(x,\omega)\in 
{\bf R}^n\times \Omega_\kappa$.
Moreover,  $Q$ satisfies the Gevrey estimates (ii) in 
${\bf R}^n\times \Omega$. 
 We set
$\psi(x,\omega)=Q(x,\omega) + \langle x, R(\omega)\rangle$. 
Recall that  $\mbox{dist}\, (\Omega_\kappa, {\bf R}^n\setminus
\Omega) \ge \kappa$. Then multiplying $Q$ and $R- \nabla g^0$ by
a suitable 
cut-off function  $h\in {\cal G}^{\rho'}_{\widetilde C}(\Omega)$, 
with $\widetilde C =
C\kappa^{-1}$ and  $C>0$
independent of $\Omega\subset \Omega_0$, such that  $h=1$ in a
neighborhood of $\Omega_\kappa$ and $h(\omega)=0$ if 
$\mbox{dist}\, (\omega, {\bf R}^n\setminus \Omega) \le \kappa/2$, we can
assume that $\psi(x,\omega) =  \langle x, \nabla g^0(\omega)\rangle$ 
for any $\omega$ such that  
$\mbox{dist}\, (\omega,{\bf R}^n\setminus \Omega) \le \kappa/2$. 
This does not change the corresponding Gevrey estimates for $\psi$. 

Let  $\varepsilon_H L_1^{N+2(\tau +2)}\le \epsilon \ll 1$. 
Then $\kappa A C_1 (C_2\kappa^{-1})L_1^{N/2}\sqrt{\varepsilon_H} \le
AC^2\epsilon \ll 1$, and  
the map 
$\Omega\ni\omega\rightarrow \nabla_x \psi(x,\omega)\in D$ becomes  a
diffeomorphism  for any $x$ fixed, which gives 
a ${\cal G}^{\rho,\rho'}$-foliation of ${\bf T}^n\times D$ 
by  Lagrangian tori 
$\Lambda_\omega\ =\ \{(p(x),\nabla_x\psi(x,\omega))\, :\, 
x\in{\bf R}^n\}\ ,\ \omega\in \Omega$.  
The action $I=(I_1,\ldots,I_n)\in D$ on each $\Lambda_\omega$,
$\omega\in D$,  is given
by
$$ 
I_j(\omega) =(2\pi)^{-1}\, \int_{\gamma_j(\omega)}\, \sigma  =
(2\pi)^{-1} (\psi(2\pi e_j,\omega) - \psi (0,\omega)) 
= R_j(\omega), 
$$ 
where $\gamma_j(\omega) = 
\{(p(te_j),\, \nabla_x\psi (te_j,\omega))\, :\, 
0\leq t\leq 2\pi\}$. 
Then $I(\omega)-\nabla g^0(\omega)= R(\omega)-\nabla g^0(\omega)$
satisfies (ii) in $\Omega$, and choosing 
$\varepsilon_H L_1^{N+2(\tau + 2)} \le \epsilon \ll 1$
we obtain that the frequency map  
$\Omega\ni\omega \mapsto I(\omega)\in D$ is a  diffeomorphism of
Gevrey class ${\cal G}^{\rho'}$. Using Remark A.1 we show  that the
inverse map $D\ni I \mapsto \omega(I)\in \Omega$ is in  ${\cal
G}^{\rho'}$ and 
$$ 
\left|\partial^\alpha \left(\omega(I) - \nabla H^0(I)\right)\right|\ \le \
\kappa   A \, 
 (C_2\kappa^{-1})^{|\alpha|}\, 
\alpha !^{\, \rho'}L_1^{N/2} \sqrt{\varepsilon_H} \, , 
$$
uniformly with respect to 
$(\varphi,\omega)\in {\bf T}^n\times \Omega$ 
and for any $\alpha,\beta \in {\bf
N}^n$.  Now we set $\Phi(x,I) = \psi(x,\omega(I))$. Then  
$\Phi(x,I) - \langle x, I\rangle$ is $2\pi$-periodic with
respect to $x$, and using
Proposition A.3     we get  the estimates 
$$ 
\left|\partial_x^\alpha \partial_I^\beta 
(\Phi(x,I) - \langle x, I\rangle)\right|\ 
\le\  \kappa   A  C_1^{|\alpha|} (C_2\kappa^{-1})^{|\beta|}  
\alpha !\, ^{\rho } \beta !\, ^{\rho(\tau + 1) +1}\, 
L_1^{N/2} \sqrt{\varepsilon_H} .
$$
Solving the equation $\Phi_I(\theta, I) = \varphi$
with respect to $\theta$ by means of Proposition  A.2  
we obtain the symplectic transformation $\chi$. For any $\omega\in
\Omega_\kappa$ and any $\theta$ we
have $(\theta, F(\theta,\omega)) = (\theta, \Phi_\theta(\theta,I(\omega)))
= \chi(\Phi_I(\theta,I(\omega)), I(\omega))$, hence, $\Lambda_\omega =
\chi({\bf T}^n\times \{I(\omega)\})$. 
Set $\tilde H(\varphi, I)=  
H(\chi(\varphi, I))$. Then $\tilde H$ is constant on 
${\bf T}^n\times \{I(\omega)\}$. We set $K(I) = \tilde H(0, I)$ and 
$R(\varphi, I) = \tilde H(\varphi, I) - K(I)$. Then $R(\varphi, I) =
0$ on ${\bf T}^n\times E_\kappa$, hence, all the derivatives of $R$
vanish on ${\bf T}^n\times E_\kappa$, since each point of $E_\kappa$
is  of positive Lebesgue density in $E_\kappa$. 
Using Proposition A.4 for $\tilde H(\varphi, I) = 
H(\theta(\varphi,I),\Phi_\theta(\theta(\varphi,I),I))$ 
we obtain for any $\alpha,\beta \in {\bf
N}^n$
$$
\left| \partial_\varphi^\alpha \partial_I^\beta 
(\tilde H(\varphi,I) - H^0(I)) 
\right|\ 
\le\  \kappa  A  C_1^{|\alpha|} (C_2\kappa^{-1})^{|\beta|}  
\alpha !\, ^{\rho } \beta !\, ^{\rho(\tau + 1) +1}\,  
L_1^{N/2} \sqrt{\varepsilon_H}
$$
uniformly with respect to 
$(\varphi,I)\in {\bf T}^n\times D$, where the constants $A,
C_1$ and $C_2$ are as above.    \finishproof

\appendix
\section*{Appendix}
\setcounter{equation}{0}

\newtheorem{theoa}{Theorem A.}
\newtheorem{lemmaa}{Lemma A.}
\newtheorem{propa}{Proposition A.}
\newtheorem{rema}{Remark A.}
\newtheorem{cora}{Corollary A.}

\renewcommand{\theequation}{A.\arabic{equation}}

We shall obtain a variant of the implicit function theorem of 
Komatsu \cite{kn:Kom} in anisotropic Gevrey classes.  
Let  $X$  and $\Omega_0$ be  domains in ${\bf R}^n$ and ${\bf R}^m$
respectevly.  
Fix $\rho' \ge \rho \ge 1$. 
Let $F=(f_1,\ldots,f_n)$ be a Gevrey function of the class 
${\cal G}^{\rho,\rho'}(X\times \Omega,  {\bf R}^n)$.  
We suppose that there are constants $A_0>0$ and $h_1, h_2 > 0$, 
and a small parameter $0<\varepsilon \le  1 $,  such that 
\begin{equation}
\forall (x,\omega)\in X\times\Omega_0 \, , \quad
\left|\partial_x^\alpha \partial_\omega^\beta 
(F(x,\omega) - x)\right| \le \varepsilon A_0 h_1 ^{|\alpha|} 
h_2 ^{|\beta|}
 \alpha !^{\, \rho}   \beta !^{\, \rho '}
                \label{A.1}
\end{equation} 
for any multi-indices $\alpha, \beta$. We choose 
$ \varepsilon >0$ so that $\varepsilon  A_0 h_1 \le 1/2$. Then   
$\left| D_x F (x,\omega) - {\rm Id} \right| \le 1/2$ for any 
$(x,\omega)\in X\times\Omega_0$, 
where $D_x F$ stands for the Jacobian of $F$ with respect to
$x$. Hence, the 
inverse matrix $(D_x F)^{-1}$ exists  and 
$ \left| D_x F(x,\omega)^{-1}\right| \le 2$ in $X\times \Omega_0$. 
Consider now a local  solution $x=g(y,\omega)$, $(y,\omega) \in
Y\times\Omega$ of $F(x,\omega)=y$, where $Y\subset {\bf R}^n$ and
$\Omega \subset \Omega_0$ are suitable domains.

\begin{propa}
Suppose that (\ref{A.1})  holds and with $0<\varepsilon \le 1$ and  
$\varepsilon A_0 h_1 \le  1/2$, and consider a local  solution $x =
g(y,\omega)$, $(y,\omega)\in Y\times \Omega$, of $F(x,\omega)=y$.
Then
$g\in {\cal G}^{\rho,\rho'}(Y\times \Omega, X)$, and there exist positive
constants $A$ and $C$ depending only on $\rho$, $\rho'$, $n$, 
and $m$, 
such that 
$$
\forall\, \alpha , \beta,  \quad 
\sup_{Y\times\Omega} \, \left|\partial_y^\alpha \partial_\omega^\beta 
(g(y,\omega) - y)\right| \le  \varepsilon A_0 A C ^{|\alpha +\beta|}\,  
h_1^{|\alpha|}\,
h_2^{|\beta|}\, \alpha !^{\, \rho }   \beta !^{\, \rho '}\, .   
$$
\end{propa}
{\em Remark A.1.} 
Let  $f=(f_1, \ldots , f_n) \in{\cal G}^{\rho}(X,{\bf R}^n)$ and 
$\left|\partial_x^\alpha 
(f(x) - x)\right| \le \varepsilon A_0 h ^{|\alpha|} 
 \alpha !^{\, \rho}$ in $X$ for any $\alpha$, 
where  $0<\varepsilon \le 1$ and  
$\varepsilon A_0 h \le  1/2$, and let $g:Y\to X$ be an inverse map to $f$.
Then the inverse map $g$ of $f$ belongs
to ${\cal G}^{\rho}(Y, X)$ and 
$\left|\partial_y^\alpha 
(g(y) - y)\right| \le  \varepsilon A_0 A C ^{|\alpha|}\,  
h^{|\alpha|}\,  \alpha !^{\, \rho }$ in $Y$ for any $\alpha$. 

\vspace{0.5cm}
\noindent  
As a corollary we obtain 
\begin{propa}
Suppose that 
$F\in {\cal G}^{\rho,\rho'}({\bf T}^n\times \Omega,{\bf T}^n )$ 
satisfies (\ref{A.1}) with   $0 <\varepsilon\le 1$ and  
$\varepsilon A_0 h_1 \le 1/2$. Then there exists  
$g \in {\cal G}^{\rho,\rho'}({\bf T}^n\times \Omega, {\bf T}^n)$ and there are 
positive constants $A$ and $C$
depending only on  $\rho$, $\rho'$, $n$, and $m$
such that
$$
\forall\, \alpha , \beta,  \quad 
\sup_{X\times\Omega} \, \left|\partial_y^\alpha \partial_\omega^\beta 
(g(y,\omega) - y)\right| \le  \varepsilon A_0 A C ^{|\alpha +\beta|} 
h_1^{|\alpha|}\, h_2^{|\beta|}\, 
\alpha !^{\, \rho }   \beta !^{\, \rho '}\, .   
$$
\end{propa}
{\em Proof of Proposition A 1}. 
We rescale the variables 
$$
x \mapsto h_1 x = \tilde x \in h_1 X := \tilde X\, , \quad 
y \mapsto h_1 y = \tilde y \in h_1 Y := \tilde Y\, , \quad 
\omega \mapsto h_2\omega = \tilde \omega \in  h_2\Omega :=
\tilde\Omega\, .
$$   
Set $\tilde F(x,\omega) = h_1 F( h_1^{-1}x,
h_2^{-1}\omega)$, $(x,\omega)\in \tilde X\times \tilde \Omega$, and 
 $\tilde g(y,\omega) = h_1 g(h_1^{-1}y,h_2^{-1}\omega)$, 
$(y,\omega)\in \tilde Y \times \tilde \Omega$. Then  $x= \tilde
g (y,\omega)$,  $(y,\omega)\in \tilde Y \times \tilde \Omega$,  is 
a solution of the equation $\tilde F(x,\omega) = y$. Moreover,
(\ref{A.1}) implies 
\begin{equation}
\sup_{\tilde X\times \tilde \Omega} \, 
\left|\partial_x^\alpha \partial_\omega^\beta 
(\tilde F(x,\omega) - x)\right| \le \tilde \varepsilon  
 \alpha !^{\, \rho}   \beta !^{\, \rho '} ,
                          \label{A.2a}
\end{equation}
where $\tilde \varepsilon := \varepsilon A_0 h_1 \le 1/2$. 
We are going to prove   
that there exist positive constants 
$A$ and $C$ depending only on $\rho$, $\rho'$, $n$, and $m$
such that 
\begin{equation}
\forall\, \alpha , \beta,  \quad 
\sup_{\tilde Y\times\tilde \Omega} \, 
\left|\partial_y^\alpha \partial_\omega^\beta 
(\tilde g(y,\omega) - y)\right| \le  \tilde \varepsilon  
A C ^{|\alpha +\beta|}\,  
 \alpha !^{\, \rho }   \beta !^{\, \rho '}\, .   
                          \label{A.2}
\end{equation}
Then rescaling back the variables  we obtain the estimates in Proposition
A.1.

We are going to prove (\ref{A.2}). From now on we omit  `$\sim$' to
simplify the notations. 
The implicit function theorem of Komatsu
\cite{kn:Kom} implies that  
$g\in {\cal G}^{\rho'}(Y\times \Omega;  X)$. 
Moreover, it
follows from \cite{kn:Kom} that 
\begin{equation}
\sup_{ X\times  \Omega} \, 
\left|\partial_{ x}^\alpha \partial_{ \omega}^\beta 
( g_j(y,\omega) - y_j)\right| \le \varepsilon A\,  h ^{|\alpha| +
|\beta|} (\alpha !\, \beta !)^{\, \rho '}\, . 
                \label{A.3}
\end{equation}
Hereafter,  $A$ and $h$ are positive constants,  
depending only on $\rho$, $\rho'$, $n$, and $m$.
To obtain (\ref{A.3}) we consider the
inverse mapping $(g,{\rm id})$ of $(F,{\rm id})$ in $X\times\Omega$
and we use that $B=2$, $0<C=\varepsilon  \le 1$ and $h=1$ in the
estimates of $b_r$ in (4),  \cite{kn:Kom}, p. 70-71 
(see also the estimates of
$\psi^{r,\alpha}$ below).  

Take any $x^0, y^0 \in X$ and $\omega^0 \in \Omega$
such that $F(x^0, \omega^0) = y^0$. Changing the
variables if necessary we assume  $y^0 = 0$ and  $
\omega^0   = 0$. 
Consider now the solution $\omega \mapsto g(0,\omega)$ of $F(x,\omega)
= 0$ with $g(0,0) =  x^0$. 
Then we obtain 
$$
F(x,\omega) = (x-g(0,\omega))D_x F(g(0,\omega),\omega) + R(x,\omega),
$$
in a neighborhood of $(0,0)$, where $R(x,\omega) =
O(|x-g(0,\omega)|^2)$, $R\in {\cal G}^{\rho,\rho'}(X\times \Omega;
{\bf R}^n)$. Moreover, it follows from (\ref{A.2a}),  (\ref{A.3}) and
from  Proposition A.3 below that 
$$
\sup_{X\times\Omega} \, \left|\partial_x^\alpha \partial_\omega^\beta 
R(x,\omega)\right| \le \varepsilon A h^{|\alpha +
\beta|} \alpha !^{\, \rho}   \beta !^{\, \rho '},
$$
where $A>0$ and $h>0$ depend only  on 
$\rho$, $\rho'$, $n$, and $n$. 
Now me make a local change of the variables 
$$
z = (x-g(0,\omega)) D_x F(g(0,\omega),\omega),
$$
and we obtain the
equation (with respect to $z$)
\begin{equation}
z = y + \varphi (z,\omega)\, ,\quad z\in Y_0,\, y\in Y,\ \omega \in \Omega\, ,
               \label{A.4}
\end{equation}
where $Y_0$ and  $Y$ are  neighborhoods of $0$ in ${\bf R}^n$,
$\varphi(0,\omega)\equiv 0$, and 
$$
\forall\, \alpha, \beta,\quad 
\sup_{Y_0\times\Omega} \, \left|\partial_z^\alpha \partial_\omega^\beta 
\varphi (z,\omega)\right| \le \varepsilon A h^{|\alpha +
\beta|} \alpha !^{\, \rho}   \beta !^{\, \rho '} 
$$
with some positive constants $A$ and $h$ as above. We are going to
estimate the derivatives $\partial_y^p\partial_\omega^\alpha u(0,0)$
of the solution $z=u(y,\omega)$ of the equation above. 

For $j\in {\bf N}$ we  
set $M_j = j !^{\, \rho}$ and $N_j = j !^{\, \rho' }$. Then 
$M_j \le N_j$ and there is $H\ge 1$ such that 
\begin{equation}
\left( M_q/q!\right)^{ 1/q-1} \le H \left( M_p/p!\right)^{ 1/p-1} ,\ 
\left( M_q/q!\right)^{ 1/q} \le H  \left( M_p/p!\right)^{ 1/p},\ 2\le
q \le p. 
                \label{A.5}
\end{equation}
The same estimates hold for $N_j$ as well. 

Consider the Taylor expansion of $\varphi$ in formal power series at
$(0,0)$ 
$$
\varphi(z,\omega) = \sum_{|p|\ge 2}\, \sum_{\alpha\in {\bf N}^m} \,
\frac{\varphi^{p,\alpha}}
{p! \alpha !} z^p \omega^\alpha\, ,\quad  \varphi^{p,\alpha} = 
\partial_z^p \partial_\omega^\alpha \varphi(0,0).
$$
Then we have
$$
\varphi \ll  \varepsilon A \sum_{|p|\ge 2}\, \sum_{\alpha\in {\bf N}^m} \,
\frac{M_{|p|}}{p!}\frac{N_{|\alpha|}}{\alpha !} 
h^{|p|+|\alpha|} z^p \omega^\alpha,
$$
which means  that 
$|\varphi^{p,\alpha}| \le  
\varepsilon A M_{|p|}N_{|\alpha|} h^{|p|+ |\alpha|}$ for each
$p,\alpha$. Consider the solution $z=u(y,\omega)$ of (\ref{A.4})  
and its formal
Taylor series at $(0,0)$ 
$$
u(y,\omega) = y + \sum_{|p|\ge 2}\, \sum_{\alpha\in {\bf N}^m}\,
\frac{u^{p,\alpha}}
{p! \alpha !} y^p \omega^\alpha. 
$$
Then (\ref{A.4}) is formally equivalent to 
$$
\sum_{|p|\ge 2}\, \sum_{\alpha\in {\bf N}^m} \,
\frac{u^{p,\alpha}}
{p! \alpha !} y^p \omega^\alpha = \sum_{|p|\ge 2}
\sum_{\alpha\in {\bf N}^m}\,
\frac{\varphi^{p,\alpha}}
{p! \alpha !} \left(y + \sum_{|q|\ge 2}\sum_{\beta\in {\bf N}^m} \,
\frac{u^{q,\beta}}
{q! \beta !} y^q \omega^\beta\right)^p \omega^\alpha . 
$$
Denote by  $v(y,\omega)$ the power  series 
$$
v(y,\omega) =  \sum_{|p|\ge 2}\sum_{\alpha\in {\bf N}^m} \,
\frac{v^{p,\alpha}}
{p! \alpha !} y^p \omega^\alpha, 
$$ 
and suppose that it is a formal solution of
\begin{equation} 
v(y,\omega) =  \varepsilon A \sum_{|p|\ge 2}\sum_{\alpha\in {\bf N}^m} \,
\frac{M_{|p|}}{p!}\frac{N_{|\alpha|}}{\alpha!}\,
h^{|p|} h^{|\alpha|}\,  (y+v(y,\omega))^p\,  \omega^\alpha \,  . 
               \label{A.6}
\end{equation}
We claim that such a solution exists and it is unique, $v^{p,\alpha}>0$, and 
$u \ll y +v$. Indeed, for any $|p|\ge 2$, $\alpha\in {\bf N}^m$, 
and $1\le j\le n$,
there is a
polynomial $Q_{p,\alpha,j}\not\equiv 0$  
of the variables $(t_{q,\beta}, s_{r,\gamma})$ 
and with non-negative coefficients, where  
$q,r\in {\bf N}^n$, $\beta, \gamma \in {\bf N}^m$,  $2\le |q|\le
|p|$, $\beta \le \alpha$, $2\le |r|\le
|p|-1$, $\gamma \le \alpha$,  and such
that $u_j^{p,\alpha}$, respectively,   $v_j^{p,\alpha}$, is the value of 
$Q_{p,\alpha,j}$ for $t_{q,\beta} =\varphi^{q,\beta} $, 
$s_{r,\gamma} = u_j^{r,\gamma}$, respectively,  for 
$t_{q,\beta} =\varepsilon A M_{|p|}N_{|\alpha|}h^{|p|+|\alpha|}$, 
$s_{r,\gamma} = v_j^{r,\gamma}$. For $|p|=2$ the polynomial
$Q_{p,\alpha,j}$ is independent of $s_{r,\gamma}$, $|r|\ge 2$. Since 
the coefficients of $Q_{p,\alpha,j}$ are non-negative, 
we get $|u_j^{r,\alpha}| \le
v_j^{r,\alpha}$ for  $|r|=2$ and each  $\alpha$. 
By recurrence we obtain 
$|u_j^{r,\alpha}| \le
v_j^{r,\alpha}$ for each $|r|\ge 2$ and  ${\alpha\in {\bf N}^m}$.

Now we set  $y=(s,\ldots,s)$, $s>0$, and consider 
$$
\phi(s,\omega) =  \sum_{r= 2}^\infty \sum_{\alpha\in {\bf N}^m} \,
\frac{\phi^{r,\alpha}}
{r! \alpha !} s^r \omega^\alpha, \ 
\phi^{r,\alpha} = \sum_{|p|= r} v^{p,\alpha}\frac{r!}{p!} . 
$$ 
Note that  for such $y$ the right hand side of the 
equation (\ref{A.6}) is invariant under any
permutation of the components $v_j(y,\omega)$, $1\le j\le n$, 
and using the uniqueness  of the solution, we obtain 
$\phi_1(s,\omega)=\cdots = \phi_n(s,\omega) = \tilde\phi(s,\omega)$. 
Then $t=s+\tilde\phi(s,\omega)$ is a formal solution of 
$$
t = s + \varepsilon A \sum_{r= 2}^\infty \sum_{\alpha\in {\bf N}^m} \,
\sum_{|p|=r} \frac{r!}{p!}\, 
\frac{M_{r}}{r!}\frac{N_{|\alpha|}}{\alpha!}
 (ht)^r (h\omega)^\alpha \, , 
$$
where $\sum_{|p|=r} \frac{r!}{p!} \le n^r $.   
Let
$t=\psi(s,\omega) := s + \sum_{r= 2}^\infty \sum_{\alpha\in {\bf N}^m} \,
\frac{\psi^{r,\alpha}}
{r! \alpha !} s^r \omega^\alpha$ be the formal solution of
\begin{equation} 
t = s + \varepsilon A  \sum_{r= 2}^\infty\,  
\sum_{\alpha\in {\bf N}^m} \,
\frac{M_{r}}{r!}\frac{N_{|\alpha|}}{\alpha!}
 (hnt)^r (h\omega)^\alpha \, . 
               \label{A.7}
\end{equation}
We obtain as above $s+ \tilde\phi \ll \psi$, which implies
$|u_j^{p,\alpha}| \le \psi^{|p|,\alpha}$ for $1\le j\le n$ and 
for any $|p|\ge 2$, $\alpha\in {\bf N}^m$.     

We are going to estimate  $\psi^{r,\alpha}$. 
Set $C(\omega) = A\,  \sum_{\alpha\in {\bf N}^m} \,
\frac{N_{|\alpha|}}{\alpha !} (h\omega)^\alpha$. Then (\ref{A.7})
becomes
$$
s =  k(t,\omega) := t - \varepsilon C(\omega) \sum_{r= 2}^\infty \,
\frac{M_{r}}{r!} (hnt))^r. 
$$
Recall that $0<\varepsilon \le 1$. We suppose also that $A\ge 1$ and
$h\ge 1$. 
As in \cite{kn:Kom}, using the Lagrange expansion theorem,  
we obtain for $r\ge 2$ 
$$
\begin{array}{rcl}
\psi^{r,\alpha} &=& \ \,  \displaystyle\left[\partial_\omega^\alpha 
\left\{\left(\frac{d}{dt}\right)^{r-1}
\left(\frac{t}{k(t,\omega)}\right)^{r}\right\}|_{t=0} \right]
|_{\omega=0} \\ [0.5cm]
&\le& 
\displaystyle\varepsilon \left[\partial_\omega^\alpha 
\left\{\left(\frac{d}{dt}\right)^{r-1}
\left( \sum_{q=0}^\infty \left(C(\omega)hn  
\sum_{p=1}^{r-1} 
\frac{M_{p+1}}{(p+1)!}(hnt)^p\right)^q \right)^{r}\right\}|_{t=0} \right]
|_{\omega=0}\,  , 
\end{array}
$$
and we get
$$
\begin{array}{rcl}
\psi^{r,\alpha}&\le& \displaystyle\varepsilon M_r \, (4Hhn)^{r-1}\,  
\left[ \partial_\omega^\alpha 
\left(1 + Ahn\sum_{\beta\le \alpha} \,
\frac{N_{|\beta|}}{|\beta| !} 
(h\omega)^\beta\right)^{r-1} \right]|_{\omega=0}\\ [0.5cm]
&\le& \displaystyle\varepsilon M_r (8AHh^2n^2)^{r-1} 
\left[\partial_\omega^\alpha 
\left(\sum_\beta \,
\left( \left(\frac{N_{|\alpha|}}{|\alpha| !}\right)^{1/|\alpha|} h\omega 
\right)^\beta\right)^{r-1}\right]|_{\omega=0}\\  [0.5cm]
&\le& \displaystyle\varepsilon M_r N_{|\alpha|} 
 (8AHh^2n^2)^{r-1} h^{|\alpha|}
2^{m(r-2)+|\alpha|} . 
\end{array}
$$
We have used the inequalities  
(\ref{A.5}) for $M_j$ and $N_j$,  as well as the
identity 
$$
\left(\sum_{\beta\in {\bf N}^m} \omega^\beta\right)^{q} = 
\sum_{\alpha\in {\bf N}^m}\, C_\alpha \omega^\alpha,
$$
where $q=r-1\ge 1$ and 
$$
 C_\alpha =  \left(\begin{array}{ccc}
                 q+\alpha_1 -1 \\ \alpha_1
                     \end{array}\right) \cdots 
\left(\begin{array}{ccc}
                 q+\alpha_m -1 \\ \alpha_m 
                     \end{array}\right) 
\ \le \ 2^{m(q-1)+|\alpha|} . 
$$
Thus we obtain 
$$
\forall\,  (p ,\alpha,j), \quad 
|u_j^{p,\alpha}|\le  \psi^{|p|,\alpha} \le \varepsilon A\, M_{|p|}
N_{|\alpha|} h^{r+|\alpha|}
$$
where $A, h>0$ depend  only on   
$\rho$, $\rho'$, $n$, and $m$. 
Changing back the variables and rescaling back 
$\Omega$, we obtain the desired estimates for $g$. 

Finally we recall the Gevrey estimates 
for the composition of two functions in anisotropic Gevrey classes.   
\begin{propa} Let 
$g \in {\cal G}^{\mu}_{C_1}(\Omega,Y)$  and 
$f\in  {\cal G}^{\rho,\mu} _{B,C_2}(X\times Y)$,
 where $\mu \ge \rho \ge 1$, $B$, $C_1$ and $C_2$ are positive
constants and $X$, $Y$ and $\Omega$ are open sets in 
${\bf R}^n$. 
Suppose that 
$\|g\|_{C_1} = A_1$ and  $\|f\|_{B,C_2} =
A_2$ in the corresponding  Gevrey norms with $A_1,\, A_2 >0$.  
Then the composition $(x,\omega) \mapsto F(x,\omega) : =
f(x,g(\omega))$ belongs to   ${\cal G}^{\rho,\mu}_{B,C}(X\times
\Omega)$, 
where 
$C = 2^{n+\mu} n^\mu  C_1 \max (1, A_1 C_2)$
and $\|F\|_{C}\le A_2$. 
\end{propa}
{\em Proof}. 
Using   the Faa de Bruno formula we write 
\[
\partial_x^\gamma\partial_\omega^\alpha F(x,\omega) =  
 \displaystyle{\sum_{1\le |\beta|=p \le |\alpha|} \ 
\frac{(\partial_x^\gamma\partial_y ^\beta f)(x,g(\omega))}{\beta !}\
\sum_{\scriptstyle{
{\alpha^1+\cdots+\alpha^p=\alpha}
\atop{|\alpha^1|\ge 1,\ldots, |\alpha^p| \ge 1}}}}\ 
\frac{\alpha !}{\alpha^1 ! \cdots \alpha^p !}\ (\partial_\omega^{\alpha^1}
g)(\omega) \cdots (\partial_\omega^{\alpha^p}g)(\omega) .  
\]
Since 
$|\partial_\omega^{\alpha}g(\omega)| \le A_1 C_1^{|\alpha|} \alpha!^{\,
\mu}$, and $j!\le 2^j (j-1)!$ for $j\ge 1$,      we estimate above 
$\left|\partial_x^\gamma\partial_\omega^\alpha F(x,\omega)\right|$ by 
\[ 
  A_2 B^{|\gamma|} (2^{\mu-1}C_1)^{|\alpha|} \gamma !\, ^\rho \, 
\alpha  !\, 
 \displaystyle{\sum_{1\le |\beta|=p \le |\alpha|} \ 
\sum_{\scriptstyle{
{\alpha^1+\cdots+\alpha^p=\alpha}
\atop{|\alpha^1|\ge 1,\ldots, |\alpha^p| \ge 1}}}}\ 
(A_1 C_2)^{p}\
\left( (|\alpha^1|-1) ! \cdots (|\alpha^p|-1) ! p !\right)
^{\, \mu-1}\  
\]
We have  $(|\alpha^1|-1) ! \cdots (|\alpha^p|-1) ! p ! \le
|\alpha| ! \le n^{|\alpha|}\alpha !$. Observe that for any $r,p\in {\bf
N}$, $r,p \ge 1$,  
the number of multi-indices 
$\gamma =(\gamma_1,\ldots, \gamma_p)\in {\bf N}^p$ 
with $|\gamma|=r$ is given by 
$\pmatrix{r + p -1\cr p-1\cr}$   
(see \cite{kn:Rod}, Sect. 1.2). 
Then  the
number of partitions $\alpha^1 +\cdots + \alpha^p = \alpha$, $\alpha^j
\in {\bf N}^n$, of 
$\alpha=(\alpha_1,\ldots, \alpha_n)$  is given  by
\[
\pmatrix{\alpha_1 + p -1\cr p-1\cr}\cdots \pmatrix{\alpha_{n} + p -1\cr
p-1\cr}\ \le \ 2^{|\alpha|} 2^{n(p-1)}\, ,
\]
which implies
\[
\left|\partial_x^\gamma \partial_\omega^\alpha F(x)\right| \le 
A_2 2^{-n} B^{|\gamma|}  
\left(2^{n+\mu}n^{\mu}C_1 \max (1,A_1C_2)\right)^{|\alpha|}\,  
\gamma !\, ^\rho \, \alpha!^{\,
\mu}
\sum_{\beta\in{\bf N}^n} 2^{-n|\beta| } 
<     A_2   B^{|\gamma|} C^{|\alpha|} \gamma !\, ^\rho \,  \alpha!^{\,
\mu}.
\] 
In the same way we obtain
\begin{propa} Let 
$g \in {\cal G}^{\rho,\mu}_{B_1,C_1}(X\times\Omega,Y)$  and 
$f\in  {\cal G}^{\rho,\mu} _{B_2,C_2}(Y\times \Omega)$,
 where $\mu \ge \rho \ge 1$, $B_j$ and  $C_j$, $j=1,2$,  are positive
constants, and $X$, $Y$,  and $\Omega$ are open sets in 
${\bf R}^n$. 
Suppose that 
$\|g\|_{B_1, C_1} = A_1$ and  $\|f\|_{B_2,C_2} =
A_2$ in the corresponding  Gevrey norms with $A_1,\, A_2 >0$.  
Then the composition $(x,\omega) \mapsto F(x,\omega) : =
f(g(x,\omega),\omega)$ belongs to   ${\cal G}^{\rho,\mu}_{B,C}(X\times
\Omega)$, 
where $B=2^{n+\rho}(2n)^{\rho} B_1 \max(1,A_1B_2)$, 
$C= C_2 + 2^{n+\rho}(2n)^{\rho}C_1\max(1,A_1B_2) $,  and 
$\|F\|_{B,C}\le A_2$. 
\end{propa}
To prove the Gevrey estimates we set $z=(x,\omega)$ and using
Leibnitz formula 
and the Faa de Bruno formula we write 
$$
\begin{array}{rcl}
\partial_x^\gamma\partial_\omega^\alpha F(x,\omega) &=&  
 \displaystyle{\sum_{\alpha'\le \alpha}
\pmatrix{\alpha\cr \alpha'\cr}   
\sum_{1\le |\beta|=p \le |\gamma +\alpha'|} \ 
\frac{(\partial_y^\beta\partial_\eta ^{\alpha-\alpha'} f)(g(z),\eta)|_
{\eta=\omega}}{\beta !}} \\ [0.3cm]
&\times&  \displaystyle{\sum_{\scriptstyle{
{\delta^1+\cdots+\delta^p=\delta=(\gamma,\alpha')}
\atop{|\delta^1|\ge 1,\ldots, |\delta^p| \ge 1}}}\ 
\frac{\delta !}{\delta^1 ! \cdots \delta^p !}\ (\partial_z^{\delta^1}
g)(z) \cdots (\partial_z^{\alpha^p}g)(z)} ,  
\end{array}
$$
where 
$\delta^j =(\gamma^j,\alpha^j)$. On the other hand,  
$$
\begin{array}{lcr}
(\gamma^1 ! \cdots \gamma^p !)^{\rho-1}(\alpha^1 ! \cdots \alpha^p!)^{\mu-1} 
\beta ! \, ^{\rho -1} \\
\le 2^{(\rho - 1)|\gamma +\alpha'|}\, 
((|\gamma^1 +\alpha^1|-1) ! \cdots (|\gamma^p + \alpha^p|-1) ! p
!)^{\rho-1}(\alpha^1 ! \cdots \alpha^p!)^{\mu-\rho} \\
 \le 2^{(\rho-1)|\gamma +\alpha'|} (2n)^{(\rho-1)|\gamma +\alpha'|} 
\gamma^\rho \alpha' !^{\mu} 
\end{array}
$$
and as above we complete the proof of the proposition. 

\vspace{0.5cm}
\noindent
{\em Proof of Lemma 3.4}. 
 First, using the
Cauchy formula, we obtain as in \cite{kn:Poe1}, Lemma A.3, 
that $f: O_{2(\tilde \upsilon + \upsilon)h}
\to {\bf C}$ is one-to-one (injective). Note that 
$2(\tilde \upsilon + \upsilon) = 1 - 4\upsilon$. To show that 
$ O_{\tilde \upsilon h} \subset f(O_{(1 -4\upsilon)h})$, we take 
$w\in O_{\tilde \upsilon h}$ and $\omega \in \Omega_\kappa$ such that 
$|w-\omega|<\tilde \upsilon h$. Set $f = {\rm id\, } -F$. Then 
$|F|_h \le \upsilon h$ and  $|DF|_{(1 -4\upsilon)h} \le 1/4$ by the Cauchy
estimates. Put $u_0 = \omega$ and $u_{k+1} = F(u_k) + w$ for $k\ge
0$. By recurrence we prove  
$$
|u_{k+1}- u_k| \le 2^{-1}4^{-k}(1-4\upsilon)h\ ,\quad 
|u_{k}- \omega | \le (2/3)(1-4^{-k})(1-4\upsilon)h .
$$ 
Taking the
limit we find $u\in  B_{r}(\omega)$, $r=(1-4\upsilon)h$, such that 
$u=F(u) + w$. The corresponding estimates of $\phi$ follow from the
arguments in \cite{kn:Poe1}, Lemma A.3. \finishproof

\end{document}